\documentclass[UTF8, reqno, 11pt, twoside]{amsart}
\usepackage{enumerate,amsmath,amssymb,mathrsfs,amsthm,amsfonts,mathtools}
\usepackage{bbm}
\usepackage{cases}
\usepackage{graphicx}
\usepackage{stmaryrd}
\usepackage{esint}
\usepackage{color}
\usepackage{mathptmx} 
\DeclareMathAlphabet{\mathcal}{OMS}{cmsy}{m}{n}
\DeclareSymbolFont{largesymbols}{OMX}{cmex}{m}{n}

\usepackage[
colorlinks=true, linkcolor=blue, citecolor=magenta]{hyperref}

\usepackage{fancyhdr}
{}

\numberwithin{equation}{section}

\newcounter{counterConstant}
\newcommand{\const}[1]{
    \addtocounter{counterConstant}{1}
    \edef#1{\arabic{counterConstant}}
}

\newcommand{\be}{\begin{eqnarray}}
\newcommand{\ee}{\end{eqnarray}}
\newcommand{\ce}{\begin{eqnarray*}}
\newcommand{\de}{\end{eqnarray*}}
\newtheorem{theorem}{Theorem}[section]
\newtheorem{lemma}[theorem]{Lemma}
\newtheorem{remark}[theorem]{Remark}
\newtheorem{definition}[theorem]{Definition}
\newtheorem{proposition}[theorem]{Proposition}
\newtheorem{example}{Example}
\newtheorem{corollary}[theorem]{Corollary}
\newtheorem{assumption}{Assumption}

\def\bt{\begin{theorem}}
\def\et{\end{theorem}}
\def\bl{\begin{lemma}}
\def\el{\end{lemma}}
\def\br{\begin{remark}}
\def\er{\end{remark}}
\def\bex{\begin{example}}
\def\eex{\end{example}}
\def\bd{\begin{definition}}
\def\ed{\end{definition}}
\def\bp{\begin{proposition}}
\def\ep{\end{proposition}}
\def\bc{\begin{corollary}}
\def\ec{\end{corollary}}
\def\ba{\begin{assumption}}
\def\ea{\end{assumption}}
\def\bpf{\begin{proof}}
\def\epf{\end{proof}}

\def\cB{{\mathcal B}}
\def\cC{{\mathcal C}}

\def\cK{{\mathcal K}}

\def\mE{{\mathbb E}}

\def\mH{{\mathbb H}}

\def\mL{{\mathbb L}}

\def\mN{{\mathbb N}}

\def\mP{{\mathbb P}}

\def\mR{{\mathbb R}}

\def\R{\mathbb R}
\def\N{\mathbb N}

\def\bP{{\mathbf P}}

\def\bP{{\mathbf P}}

\def\bE{{\mathbf E}}
\def\1{{\mathbf{1}}}

\def\sF{{\mathscr F}}

\def\sP{{\mathscr P}}

\def\tL{\widetilde{\mathbb{L}}}
\def\tH{\widetilde{\mathbb{H}}}

\def\div{\mathord{{\rm div}}}
\def\eps{\varepsilon}
\def\d{\text{\rm{d}}}

\def\e{\mathrm{e}}

\def\a{\alpha}
\def\om{\omega}

\def\p{\partial}

\def\l{\lambda}
\def\si{\sigma}

\def\({{\Big(}}
\def\){{\Big)}}
\def\[{{\Big[}}
\def\]{{\Big]}}
\def\<{{\langle}}
\def\>{{\rangle}}

\def\osc{\mathop{{\rm osc}}}

\def\dif{{\mathord{{\rm d}}}}

\def\l{\left}
\def\r{\right}
\def\geq{\geqslant}
\def\leq{\leqslant}

\allowdisplaybreaks

\begin{document}

\title{SDEs with critical time dependent drifts: weak solutions}

\author{Michael R\"ockner and Guohuan Zhao
} 
 
 \address{Michael R\"ockner: Department of Mathematics, Bielefeld University, Germany \\
and Academy of Mathematics and Systems Science, Chinese Academy of Sciences (CAS),
Beijing, 100190, P.R.China\\  
Email: roeckner@math. uni-bielefeld.de
 }
 
\address{Guohuan Zhao: Department of Mathematics, Bielefeld University, Germany \\ 
Email: zhaoguohuan@gmail.com
 }

\thanks{Research of Michael and Guohuan is supported by the German Research Foundation (DFG) through the Collaborative Research Centre (CRC) 1283 Taming uncertainty and profiting from randomness and low regularity in analysis, stochastics and their applications.}

\begin{abstract} 
We prove the unique weak solvability of time-inhomogeneous stochastic differential equations with additive noises and drifts in critical Lebesgue space $L^{q}([0,T]; L^{p}(\mathbb{R}^d))$ with $d/p+2/q=1$. The weak uniqueness is obtained by solving corresponding Kolmogorov's backward equations in some second-order Sobolev spaces, which is analytically interesting in itself.
\end{abstract}

\maketitle



\section{Introduction}
The main aim of this paper is to investigate the well-posedness of the following stochastic differential equation (SDE): 
\be\label{Eq-SDE}
\d X_{t} = b(t, X_{t})\d t+\sqrt{2}\d W_t,\quad X_{0}=x\in \R^d, 
\ee
where $W$ is a $d$-dimensional standard Brownian motion and $b: [0,T]\times \R^d\to \R^d$ is a vector field in some critical Lebesgue spaces $L^{q}([0,T]; L^{p}(\R^d))$ with $d/p+2/q=1$. 

The study of classical strong solutions to SDEs in multidimensional spaces with singular drifts at least date back to \cite{veretennikov1980strong2}, where Veretennikov showed that \eqref{Eq-SDE} admits a unique strong solution, provided that $b$ is bounded measurable. Using Girsanov's transformation and results from PDEs, Krylov-R\"ockner \cite{krylov2005strong} obtained the existence and uniqueness of strong solutions to \eqref{Eq-SDE}, when $b$ satisfies the following subcritical Ladyzhenskaya-Prodi-Serrin (LPS) type condition:
\be\label{Con-sub}
b \in \mL^p_q(T):=L^q([0,T]; L^p(\R^d))  \ \mbox{ with } p,q \in (2, \infty), \ \  \frac{d}{p}+\frac{2}{q}<1. 
\ee
After that, various works were devoted to generalize the well-posedness results and study the properties of solutions to SDEs with singular coefficients, among which we quote  \cite{fedrizzi2011pathwise}, \cite{lee2017existence}, \cite{menoukeu2013variational}, \cite{mohammed2015sobolev}, \cite{rezakhanlou2014regular}, \cite{xia2020lqlp}, \cite{zhang2005strong}, \cite{zhang2011stochastic}. 

However, for the critical regime:
\be\label{Con-critical}
b \in \mL^p_q(T)  \ \mbox{ with } p,q \in [2, \infty], \ \  \frac{d}{p}+\frac{2}{q}=1, 
\ee 
it has been a long-standing problem whether  SDE \eqref{Eq-SDE} is well-posed or not in the strong sense under the critical condition \eqref{Con-critical}. Beck-Flandoli-Gubinelli-Maurelli \cite{beck2019stochastic} showed that if $b$ satisfies \eqref{Con-critical} with $p>d$ or $p=d$ and $\|b\|_{\mL^d_\infty(T)}$ is sufficiently small, \eqref{Eq-SDE} has at least one strong solution starting from a diffusive random variable in a certain class. In \cite{nam2020stochastic}, Nam proved that \eqref{Eq-SDE} admits a unique strong solution for each $x\in \R^d$ when the Lebesgue-type $L^q$ integrability in the time variable is replaced by a stronger Lorentz-type $L^{q,1}$ integrability condition:
$$
b \in L^{q,1}([0,T]; L^p(\R^d))  \ \mbox{ with } p,q \in (2, \infty), \ \  \frac{d}{p}+\frac{2}{q}=1. 
$$
But the endpoint cases $(p,q)=(d,\infty)$ and $(p,q)=(\infty,2)$ are excluded in \cite{nam2020stochastic}. Very recently, Krylov made significant progress in his works \cite{krylov2020strong}, where the strong well-posedness is proved in the case that $b(t,x)=b(x)\in L^d(\R^d)$ with $d\geq 3$. His approach is based on his earlier work with Veretennikov \cite{veretennikov1976explicit} about the Wiener chaos expansion for strong solutions of \eqref{Eq-SDE}, and also some new estimates presented in \cite{krylov2020stochastic2} and  \cite{ krylov2020stochastic1}. 

Surprisingly, to the authors' best knowledge, there is no complete answer even for the unique {\em weak} solvability of \eqref{Eq-SDE} at the Lebesgue-critical regime \eqref{Con-critical}. There are few works on this subject, among them, we mention that  Wei-Lv-Wu  \cite{wei2017weak} studied the weak well-posedness of \eqref{Eq-SDE} when $b\in \cC_{q}^0([0,T]; L^p(\R^d))$ with $p,q<\infty$. Recently, Kinzebulatov-Sem{\"e}nov \cite{kinzebulatov2019brownian} not only proved the weak existence of solutions to \eqref{Eq-SDE}, but also constructed a corresponding Feller process when $b(t,x)=b(x)$ is form-bounded, which includes  the case that $b$ is in the weak $L^d$ space and $\|b\|_{L^{d,\infty}}$ is small enough (see also \cite{kinzebulatov2020feller}). Almost simultaneously, 
Xia-Xie-Zhang-Zhao \cite{xia2020lqlp} proved the weak well-posedness of \eqref{Eq-SDE} when $b\in C([0,T]; L^{d}(\R^d))$. However, the borderline case $b\in \mL^d_\infty(T)$ is much more delicate and cannot be solved by the arguments used in \cite{kinzebulatov2019brownian} or \cite{xia2020lqlp}. In \cite{zhang2020stochastic}, Zhang and the second named author of this paper studied \eqref{Eq-SDE} at the supercritical regime, and they proved that if $b,\, \div b \in \mL^p_q(T)  \ \mbox{ with } p,q \in [2, \infty)\ \mbox{ and }\   \frac{d}{p}+\frac{2}{q}<2$, then SDE \eqref{Eq-SDE} has at least one weak (martingale) solution.  We also mention that in a series of very recent works \cite{krylov2020diffusion}, \cite{krylov2020stochastic2},  \cite{krylov2020time1} and \cite{krylov2020time2}, when $b\in L^{d+1}(\R^{d+1})$ or $b\in L^d(\R^d)$, Krylov not only constricted the strong Markov processes associated with  \eqref{Eq-SDE}, but also studied many properties of these processes, such as Harnack's inequality, higher summability of Green's functions, and so on.

Now let us state the main motivations of this paper. The significance of the LPS condition \eqref{Con-sub} and \eqref{Con-critical} named after the authors who posed these conditions to prove global well-posedness of $3$D Navier-Stokes equations and smoothness of solutions.  For the regularity theory of Navier-Stokes equations, the endpoint case $(p, q) = (3,\infty)$, which triggered a lot of papers, is extremely difficult,  and was finally solved by Escauriaza-Seregin-{\v{S}}ver{\'a}k in \cite{escauriaza2003backward} (see also \cite{dong2009navier} and \cite{gallagher2013profile}). As presented in \cite{zhang2016stochastic} and  \cite{zhao2019stochastic}, by letting $b$ in \eqref{Eq-SDE} be a solution to the Navier-Stokes equation (which is a divergence free vector field), the stochastic equation \eqref{Eq-SDE} can be related with the Navier-Stokes equation through  Constantin and Iyer's representation (see \cite{constantin2008stochastic} and \cite{zhang2010astochastic}). This deep connection between singular SDEs and Navier-Stokes equations is our first  motivation to study \eqref{Eq-SDE} under the critical condition \eqref{Con-critical}, especially for the borderline case $(p, q)=(d,  \infty)$. Our work is also motivated by the following interesting example:  Let
$$
b(x)=- \lambda x/|x|^2, \quad x\in \R^3 \mbox{ and } \lambda>0. 
$$ 
Obviously, $b\notin L^3_{loc}(\R^3)$ but $b\in L^{3,\infty}(\R^3)$ ($=$ weak $L^3(\R^3)$ space). It was discussed in \cite{beck2019stochastic} and  \cite{zhang2020stochastic} that \eqref{Eq-SDE} has no weak solution if $\lambda$ is  large. However, when $\lambda$ is sufficiently small, it was shown by Kinzebulatov-Sem\"enov in \cite{kinzebulatov2019brownian} (see also \cite{kinzebulatov2020feller}) that \eqref{Eq-SDE} has at least one weak solution. In this paper, we will also give the weak uniqueness for this example. 

\medskip
Denote the localized $\mL^p_q(T)$ (weak $\mL^p_q(T)$) space by $\widetilde \mL^p_q(T)$ $(\widetilde \mL^{p,\infty}_q(T))$ (see Section \ref{Sec-Pre}  for the precise  definitions). Our main result is 
\bt\label{Th-Main}
Assume $d\geq 3$ and that $b$ satisfies one of the following two assumptions: 
\begin{enumerate}[(a)]
\item $b=b_0+b_1$,  where  $b_1\in \mL^{p_1}_{q_1}(T)$ with $d/p_1+2/q_1=1$ and $p_1\in (d, \infty)$, and $b_0\in \widetilde \mL^{d,\infty}_\infty(T)$ with $\|b_0\|_{\widetilde \mL^{d,\infty}_\infty(T)}\leq \eps$, for some constant $\eps>0$ only depending on $d, p_1, q_1$; 
\item $b\in \tL^{d,\infty}_\infty(T)$ and $\div b\in \tL^{p_2}_{\infty}(T)$ with $p_2\in (d/2,\infty)$. 
\end{enumerate}
Then there is a unique weak solution to \eqref{Eq-SDE} such that the following Krylov type estimate is valid: 
\be\label{Eq-Krylov}
\bE \l(\int_{0}^{T} f(t, X_t) \d t \r) \leq C \|f\|_{\widetilde \mL^{p_3}_{q_3}(T)}, \quad \mbox{ for any } p_3,q_3\in (1,\infty) \mbox{ with } \frac{d}{p_3}+\frac{2}{q_3}<2. 
\ee
Here $C$ is a constant, which does not depend on $f$. 
\et

\br
\
\begin{enumerate}
\item If $b\in C([0,T]; L^d)$, then for any $\eps>0$, there exist two functions  $b_0$ and $b_1$ such that $\|b_0\|_{\mL^d_\infty(T)}<\eps$ and $b_1\in \mL^d_\infty(T)\cap L^\infty([0,T]\times \R^d)$. Therefore, any functions in $C([0,T]; L^d)$ satisfies condition (a)  in Theorem \ref{Th-Main} for arbitrary $\eps>0$. 
\item As shown in \cite[Theorem 5.1]{zhao2019stochastic}, for any $p\in (d/2, d)$, there exists a divergence free vector field $b\in L^p+L^\infty$ such that weak uniqueness of \eqref{Eq-SDE} fails. So our condition (b) is optimal when $b$ is divergence free. 
\item  Our result also can be extended to SDEs driven by multiplicative noises: 
$$
\d X_{t} = b(t, X_{t})\d t+\sigma(t, X_t)\d W_t,\quad X_{0}=x\in \R^d, 
$$ 
provided that $a:=\frac{1}{2}\si\si^t$ is uniformly elliptic and uniformly continuous in $x$ with respect to $t$, and $b, \p_j a^{ij}$ meet the same conditions as $b$ in Theorem \ref{Th-Main}. 

\item Let $d=3$, $\lambda\in \R$, 
$$
b(x)=\l(\frac{\lambda x_1x_3}{(x_1^2+x_2^2)|x|}, \frac{\lambda x_2x_3}{(x_1^2+x_2^2)|x|}, \frac{-\lambda}{|x|}\r). 
$$
We note that in this case $b\in L^{3,\infty}(\R^3)$ and $\div b\equiv0$. So our result implies that equation \eqref{Eq-SDE} has a unique weak solution for any $\lambda\in \R$. This and the example before Theorem \ref{Th-Main} ($b(x)=-\lambda x/|x|^2$) show that the integrability (or singularity) of the drift is not the only discrimination for the well-posedness of \eqref{Eq-SDE}. The structure of the vector field also affects it. 
\end{enumerate}
\er

\medskip

Now let us explain the approach used in this paper. Take $T>0$, $f\in C_c^\infty(\R^{d+1})$ and consider the equation
\be\label{Eq-PDE}
\p_t u-\Delta u-b\cdot\nabla u=f  \mbox{ in } (0,T)\times \R^d, \quad u(0)=0.  
\ee
The existence of weak solutions to \eqref{Eq-SDE} follows from a standard tightness argument and a global maximum principle for weak solutions to \eqref{Eq-PDE}. For uniqueness, under the same conditions as in Theorem \ref{Th-Main}, we shall find a  uniformly bounded and sufficiently regular  solution $u$ to the above parabolic equation so that a generalized It\^o formula can be applied to any weak solution of \eqref{Eq-SDE} and the function $u(T-t,x)$ to obtain 
$$
u(0, X_T)-u(T, x)= - \int_0^T f(t, X_t) \d t + \sqrt{2}\int_0^T \nabla u(t, X_t)  \cdot \d W_t.
$$
Then by taking expectations of both sides, we obtain 
$$
\bE \int_0^T f(t, X_t) \d t= u(T, x),
$$ which is enough to guarantee the uniqueness of $X$ in law. Under condition (a) in Theorem \ref{Th-Main}, the solvability of \eqref{Eq-PDE}  in some second-order Sobolev spaces are proved in Theorem \ref{Th-PDEmain1} by a perturbation argument, together with a parabolic type Sobolev inequality. Like the regularity theory for $3$D Navier-Stokes equations, the endpoint case $b\in \mL^d_\infty(T)$ (without the smallness  condition on $\|b\|_{\mL^d_\infty(T)}$) is more delicate and we have no answer to the full borderline case without an additional assumption on $b$. However, when $\div b\in \tL^{\frac{d}{2}+\delta}_\infty\, (\delta>0)$, by means of De Giorgi's method, we can show that any bounded weak solution of \eqref{Eq-PDE} is indeed H\"older continuous. After that, we use another interpolation inequality of Nirenberg \eqref{Eq-GNIeq2} involving H\"older norms to show that 
$b\cdot\nabla u\in \widetilde \mL^{p_3}_{q_3}(T)$ with some $p_3,q_3\in (1,\infty)$ and $d/p_3+2/q_3<2$.  This yields that the bounded weak  solution $u$ to \eqref{Eq-PDE} is indeed in $\widetilde \mH^{2,p_3}_{q_3}(T)$ (see Theorem \ref{Th-PDEmain2}), which is regular enough to apply the generalized It\^o's formula. The above mentioned analytic results seem also to be new, and are thus of independent interest. 

\medskip

We close this section by emphasizing  again that whether \eqref{Eq-SDE} admits a unique strong solution under critical LPS condition \eqref{Con-critical} is a challenging question. One of the obstacles is that Zvonkin's type of changing  variables (cf.  \cite{zvonkin1974transformation}), which works very well for the subcritical case seems very hard to be applied under condition  \eqref{Con-critical} and might be even not possible. A possible way to overcome this is using similar arguments as in \cite{krylov2020strong}, but one needs to overcome many difficulties due to the fact that we are in the time-inhomogeneous case. In our forthcoming work \cite{rockner2021critical}, for the non-endpoint case, based on a compactness criterion for $L^2$ random fields in Wiener spaces and new estimates for some functionals of the solutions to \eqref{Eq-SDE}, we construct strong solutions directly without using Yamada-Watanabe principle. Therefore, due to a fundamental result of Cherny \cite{cherny2002uniqueness}, our weak uniqueness results in this paper will also play a role in the strong well-posedness of \eqref{Eq-SDE} with drift term in critical Lebesgue spaces. 

\section{Preliminary}\label{Sec-Pre}
In this section, we introduce some notations and present some lemmas, which will be frequently used in this paper. 

Let $D$ be an open subset of $\R^d$. For any $p\in [1,\infty]$, by $L^{p,\infty}(D)$ we mean the weak $L^p(D)$ space with finite quasi-norm given by 
$$
\|f\|_{L^{p,\infty}(D)}:= \sup_{\lambda>0}\lambda \big|\{x\in D: |f(x)|>\lambda\}\big|^{1/p}. 
$$  
For any $r>0$, we define 
$$
Q_r(t,x)= (t-r^2, t)\times B_r(x), \quad  Q_r=Q_r(0,0). 
$$
Let $I$ be an open interval in $\R$ and $Q=I\times D$. For any $p, q\in [1,\infty]$, by $\mL^p_q(Q)$  and $\mL^{p,\infty}_q(Q)$ we mean the space of functions on $Q$ with finite norm given by 
$$
\|u\|_{\mL^{p}_q(Q)}:= \|\|u(t,\cdot)\|_{L^p(D)} \|_{L^q(I)}\ \mbox{ and } \ \|u\|_{\mL^{p,\infty}_q(Q)}:= \|\|u(t,\cdot)\|_{L^{p,\infty}(D)} \|_{L^q(I)}. 
 $$
We write $u\in V(Q)$ if  
$$\|u\|_{V(Q)}^2:= \|u\|_{\mL^2_\infty(Q)}^2+ \|\nabla u\|_{\mL^2_2(Q)}^2<\infty, 
$$
and $u\in V^{0}(Q)$ if $u\in V(Q)$  and for any $t\in I$, 
$$
\lim_{h\to 0} \|u(t+h,\cdot)-u(t,\cdot)\|_{L^2(D)} =0. 
$$
Given a constant  $T>0$, {\em with a little abuse of notations}, for each $p,q\in[1,\infty]$, we set  
$$
\mL^p_q(T):=L^q([0,T];L^p(\mR^d))\ \mbox{ and }\  \mL^p(T):= \mL^p_p(T). 
$$
For $p,q\in (1,\infty), s\in \mR$, define 
$$
\mH^{s, p}_q(T)=L^q([0,T]; H^{s, p}(\mR^d)), 
$$
where $H^{s, p}$ is the Bessel potential space. The usual energy space is defined as the following way: 
$$
V(T):=\Big\{f\in \mL^2_\infty(T)\cap L^2([0,T]; H^1): \|f\|_{V(T)}:=\|f\|_{\mL^2_\infty(T)}+\|\nabla_xf\|_{\mL^2(T)}<\infty\Big\}.
$$
In this paper, we will also use the localized versions of the above functional spaces. Throughout this paper we fix a cutoff function
\be\label{Eq-chi}
\chi\in C^\infty_c(\mR^{d};[0,1]) \mbox{ with $\chi|_{B_1}=1$ and $\chi|_{B^c_2}=0$}. 
\ee
For $r>0$ and $x\in\mR^{d}$, let $\chi^{y}_r(x):=\chi\l(\frac{x-y}{r}\r)$. For any $p,q\in [1,\infty]$, define 
\begin{align*}
\widetilde \mL^p_q(T):= \l\{f\in L^q([0,T]; L^p_{loc}(\mR^d)): \|f\|_{\widetilde \mL^p_q(T)}:= \sup_{y\in \mR^d} \| f\chi^y_1 \|_{\mL^p_q(T)}<\infty.  \r\}
\end{align*}
and set $\mL^p(T):= \widetilde \mL^p_p(T)$. Similarly, 
$$
\widetilde \mL^{p,\infty}_q(T):= \l\{f\in L^q([0,T]; L^{p,\infty}_{loc}(\mR^d)): \|f\|_{\widetilde \mL^{p,\infty}_q(T)}:= \sup_{y\in \mR^d} \| f\chi^y_1 \|_{\mL^{p,\infty}_q(T)}<\infty.  \r\}. 
$$
The localized Bessel potential spaces and  energy spaces are defined as following: 
\begin{align*}
\widetilde \mH^{s,p}_q(T):= \left\{f\in L^q([0,T]; H^{s, p}_{loc}): \|f\|_{\widetilde \mH^{s,p}_q(T)}:= \sup_{y\in \R^d}\|f\chi^y_1\|_{\mH^{s,p}_q(T)} \right\}, 
\end{align*}
\begin{align*}
\widetilde V (T):=\left\{f\in \widetilde \mL^2_\infty(T)\cap\widetilde \mH^{1,2}_{2}(T):
 \| f\|_{\widetilde V (T)}:=\| f\|_{\widetilde \mL^2_\infty(T)}+\|\nabla_xf\|_{\widetilde \mL^2(T)}<\infty\right\}, 
 \end{align*}
 \begin{align*}
 \widetilde V ^0(T):=\Big\{f\in \widetilde V (T):   \mbox{for any }  y&\in \mR^d, t\mapsto f(t)\chi^y_1 \\
 & \mbox{is continuous from $[0,T]$ to $L^2(\mR^d)$} \Big\}. 
\end{align*} 

\medskip
The following  Nirenberg's interpolation inequality involving H\"older norms and De Giorgi's isoperimetric inequality are well-know (cf.  \cite{nirenberg1959extended}, \cite[Theorem 2]{kufner1995interpolation} and  \cite[Lemma 1.4]{caffarelli2010giorgi}). 
\bl[Nirenberg's interpolation inequalities]
Suppose $d\geq 2$, $j, m\in \mN$, $0<j<m$, and $1<p<q<\infty$ such that 
$$
\frac{pm}{j}<q\leq \frac{p(m-1)}{j-1}, \  \a= \frac{jq-mp}{q-p} \ \mbox{ and } \ \theta=\frac{p}{q}. 
$$
Then, 
\be\label{Eq-GNIeq2}
\|\nabla^j u\|_q \leq C \|\nabla^m u\|_p^\theta \cdot [u]_\a^{1-\theta}, \quad \forall u\in H^{m, p}\cap C^\a. 
\ee
\el
\bl[De Giorgi's isoperimetric inequality]\label{Le-De Giorgi}
There exists a constant $c_d>0$ depending only on $d$ such that the following holds. For any function $u: \R^d\to \R$, set 
$$
A_u:=\{u\geq 1/2\}\cap B_1, \quad B_u:=\{u\leq 0\}\cap B_1, \quad D_u:=\{u\in (0,1/2)\}\cap B_1, 
$$
then 
$$
\|\nabla u^+\|_2^2 \geq c_d \frac{|A_u|^2|B_u|^{2-\frac{2}{d}}}{|D_u|}. 
$$
\el
The following conclusion is a variants of Theorem 1.1 in  \cite{krylov2001heat}.  
\const{\CH2}
\bl\label{Le-heat}
Let $p,q\in (1,\infty)$, $\lambda>0$. For each $u\in L^q(\R; H^{2,p}(\R^d))\cap H^{1,q}(\R; L^{p}(\R^d))$, it holds that 
\be\label{Eq-CZest}
\|\p_t u\|_{\mL^{p}_q} + \lambda \|\nabla^2 u\|_{\mL^{p}_q} \leq C \|\p_tu-\lambda \Delta u\|_{\mL^{p}_q}, 
\ee
where $C$ only depends on $d, p, q$. 
\el
The following lemma about the $L^qL^p$-maximal regularity estimates will be used several times later. 
\bl\label{Le-heat}
Let $p,q\in (1,\infty)$, $\a\in \R$. Assume that $f\in  \widetilde \mH^{\a, p}_q(T)$, then the following heat equation admits a unique solution in $\widetilde \mH^{2+\a,p}_q(T)$: 
\be\label{Eq-Heat}
\p_t u-\Delta u= f \ \mbox{ in } (0,T)\times \R^d, \quad u(0)=0 
\ee
and
\be\label{Eq-CZest2}
\|\p_t u\|_{\widetilde \mH^{\a,p}_q(T)}+\|u\|_{\widetilde \mH^{2+\a,p}_q(T)}\leq C \|f\|_{\widetilde \mH^{\a,p}_q(T)}, 
\ee
where $C=C(d, p, q, T)$. In particular, if $f\in  \mH^{\a, p}_q(T)$,  then 
\be\label{Eq-Heat-est}
\|\p_t u\|_{\mH^{\a, p}_q(T)}+\|\nabla^2 u\|_{\mH^{\a, p}_q(T)}\leq C_{\CH2} \|f\|_{\mH^{\a, p}_q(T)}, 
\ee
where $C_{\CH2}=C_{\CH2}(d, p, q)$. 
\el
\bpf
The estimate \eqref{Eq-Heat-est} is a consequence of Theorem 1.1 in \cite{krylov2001heat}. For any $y\in \R^d$, by \eqref{Eq-Heat}, we have 
$$
\p_t (u\chi^y_1 ) -\Delta (u\chi^y_1)= f\chi^y_1 -2\nabla u\cdot\nabla\chi^y_1-u\Delta \chi^y_1.
$$
Thus,  due to \cite[Theorem 1.2]{krylov2001heat} and  \cite[Proposition 4.1]{zhang2020stochastic}, for each $t\in [0,T]$
\begin{align*}
&\|\p_t u\|_{\widetilde \mH^{\a, p}_q(t)} + \|u\|_{\widetilde \mH^{2+\a, p}_q(t)} \\
\leq& C \sup_{y\in \R^d}\l(
\|\p_t (u\chi^y_1)\|_{\mH^{\a, p}_q(t)} + \|u\chi^y_1\|_{\mH^{2+\a, p}_q(t)} \r)\\
\leq &C \sup_{y\in \R^d}\l( \|f\chi^y_1\|_{\mH^{\a, p}_q(t)}+ \|\nabla u\cdot \nabla \chi_1^y\|_{\mH^{\a, p}_q(t)}+ \|u\Delta \chi_1^y\|_{\mH^{\a, p}_q(t)} \r) \\
\leq & C  \|f\|_{\widetilde \mH^{\a, p}_q(t)}+ C \sup_{y\in \R^d}\l( \|\nabla (u\chi_2^y) \cdot \nabla \chi_1^y\|_{\mH^{\a, p}_q(t)}
+ \|(u\chi_2^y)\Delta \chi_1^y\|_{\mH^{\a, p}_q(t)} \r) \\
\leq &C \l( \|f\|_{\widetilde \mH^{\a, p}_q(t)}+ \|u\|_{\widetilde \mH^{\a+1,p}_q(t)}\r).
\end{align*}
A basic interpolation inequality yields, 
\begin{align}\label{eq-W2}
\|\p_t u\|_{\widetilde \mH^{\a, p}_q(t)} + \|u\|_{\widetilde \mH^{2+\a, p}_q(t)}\leq C_T \l(\|f\|_{\widetilde \mH^{\a, p}_q(t)}+ \|u\|_{\widetilde \mH^{\a, p}_q(t)}\r), \quad \forall t\in [0,T]. 
\end{align}
Since for any $t\in [0,T]$, 
$$
\| u(t) \chi_1^y\|_{H^{\a, p}} \leq  \int_0^t \| \p_t u(r) \chi_1^y\|_{H^{\a,p}} \, \d r\leq C_T  \l(\int_0^t \|\p_t u(r)\chi_1^y\|_{H^{\a, p}}^q \d r\r)^{1/q},
$$
together with \eqref{eq-W2}, we obtain 
$$
\|u(t)\|_{\widetilde H^{\a, p}}\leq C \|f\|_{\widetilde \mH^{\a,p}_q(T)} + C\l( \int_0^t \|u(r)\|_{\widetilde H^{\a, p}}^q \d r\r)^{1/q}. 
$$
By Gronwall's inequality, we obtain 
$$
\sup_{t\in [0,T]}\|u(t)\|_{\widetilde H^{\a, p}} \leq C \|f\|_{\widetilde \mH^{\a,p}_q(T)}, 
$$
which together with \eqref{eq-W2} yields the desired estimate. 
\epf

Next we attempt to prove a parabolic version of Sobolev inequality, which will play a crucial role in the proof of our main result. This goal can be achieved by using the Mixed Derivative Theorem, which goes back to the work of Sobolevskii's  (cf. \cite{sobolevskii1977fractional}). 

\medskip

Let $X$ be a Banach space and let $A: D(A) \to X$ be a closed, densely defined linear operator with dense range. Then $A$ is called sectorial, if 
$$
(0, \infty) \subseteq \rho(-A) \quad \text { and } \quad\left\|\lambda(\lambda+A)^{-1}\right\|_{\mathcal{L}(X)} \leq C, \quad \lambda>0,  
$$
where $\rho(-A)$ is the resolvent set of $-A$. 
Set 
$$
\Sigma_{\phi}:=\{z \in \mathbb{C} \backslash\{0\}:|\arg z|<\phi\}. 
$$
We call
$$
\phi_{A}:=\inf \left\{\phi \in[0, \pi): \Sigma_{\pi-\phi} \subseteq \rho(-A), \sup _{z \in \Sigma_{\pi-\phi}}\left\|z(z+A)^{-1}\right\|_{\mathcal{L}(X)}<\infty\right\}
$$
the spectral angle of $A$. For all $\theta\in (0,1)$, the formulas 
\be\label{Eq-Atheta}
A^\theta x= \frac{\sin \theta\pi}{\pi} \int_0^\infty \lambda^{\theta-1} (\lambda +A)^{-1} A x \ \d \lambda, \quad x\in D(A)
\ee
\be\label{Eq-A-theta}
A^{-\theta} x= \frac{\sin \theta\pi}{\pi} \int_0^\infty \lambda^{-\theta} (\lambda+A)^{-1}x \  \d \lambda, \quad x\in X
\ee
is valid (see \cite{sobolevskii1977fractional}). 
\bl[Mixed Derivative Theorem]\label{Le-MixD}
Let $A$ and $B$ are two sectorial operators  in a Banach space $X$ with spectral angels $\phi_A$ and $\phi_B$, which are commutative  and satisfy the parabolicity condition $\phi_A+\phi_B<\pi$. The coercivity estimate 
$$
\|A x\|_{X}+\lambda \|B x\|_{X} \leq M\|A x+\lambda B x\|_{X}, \quad \forall x \in D(A) \cap D(B), \, \lambda >0
$$
implies the estimate
$$
\left\|A^{(1-\theta)} B^{\theta} x\right\|_{X} \leq C\|A x+B x\|_{X}, \quad \forall x \in D(A) \cap D(B), \, \theta \in[0,1]. 
$$
\el

The above result implies the following important parabolic type Sobolev inequality. 
\const{\CPSobolevone}
\bl\label{Le-Inter}
Let $p, q\in (1,\infty)$, $r\in(p,\infty)$, $s\in (q,\infty)$. Assume $\p_t u\in \mL^{p}_q(T)$, $u\in \mH^{2,p}_q(T)$ and $u(0)=0$.  If $1<d/p+2/q=d/r+2/s+1$, then 
\be\label{Eq-PSobolev1}
\|\nabla u\|_{\mL^r_s(T)}\leq C_{\CPSobolevone} \l( \|\p_{t} u\|_{\mL_{q}^{p}(T)} + \|\nabla^2 u\|_{\mL_{q}^{p}(T)}\r), 
\ee
where $C_{\CPSobolevone}$ depends on $d, p, q,r,s$.
\el
\bpf
Let $X= L^q(\R; L^p(\R^d))$, $A=\p_t$ and $B=-\Delta$ in Lemma \ref{Le-MixD}. It is well-known that 
$$
\phi_A=\frac{\pi}{2}, \quad \phi_B=0. 
$$
Due to \eqref{Eq-CZest}, we have  
\begin{align*}
\|\p_t u\|_{\mL^p_q}+\lambda \|\Delta u\|_{\mL^p_q}\leq  C \|\p_t u-\lambda \Delta u\|_{\mL^p_q},  
\end{align*}
for all $\lambda>0$ and $C$ only depends on $d, p,q$. 
Thanks to Lemma \ref{Le-MixD}, we obtain 
\begin{align}\label{eq-frac}
\|\p_t^{1-\theta}(-\Delta)^{\theta} u\|_{\mL^p_q} \leq C \|\p_t u-\Delta u\|_{\mL^p_q}\leq C\l( \|\p_t u\|_{\mL^p_q}+ \|\nabla^2 u\|_{\mL^{p}_q} \r), 
\end{align}
for all $u\in H^{1,q}\left(\mathbb{R}, L^{p}\right) \cap L^{q}\left(\mathbb{R}, H^{2,p}\right)$.  By  \eqref{Eq-A-theta}, we have 
\be
\begin{aligned}
\p_t^{-1+\theta} f(t,x)=& \frac{\sin (1-\theta)\pi}{\pi} \int_0^\infty \lambda^{-1+\theta} (\lambda+\p_t)^{-1}f(t,x) \  \d \lambda\\
=& \frac{\sin (1-\theta)\pi}{\pi} \int_0^\infty \lambda^{-1+\theta}  \int_{-\infty}^t \e^{-\lambda (t-s)} f(s,x)\d s\  \d \lambda \\
=& \frac{\Gamma(\theta) \sin (1-\theta)\pi}{\pi}\int_{-\infty}^t (t-s)^{-\theta} f(s,x) \d s=:c_\theta (h_\theta*_t f)(t),  
\end{aligned}
\ee
where $h_\theta(t):= t^{-\theta}\1_{(0,\infty)}(t)$. 
Set 
$$
\theta=1+\frac{1}{s}-\frac{1}{q}=\frac{1}{2}+\frac{d}{2p}-\frac{d}{2r}\in \l(\frac{1}{2},1\r).
$$
Noting that $h_\theta \in L^{\frac{1}{\theta},\infty}(\R)$, by a refined version of Young's inequality (cf. \cite[Theorem 1.5]{bahouri2011fourier}),  
$$
\| \p_t^{-1+\theta} f \|_{L^s(\R)}\leq C\|h\|_{L^{1/\theta,\infty}(\R)} \|f\|_{L^q(\R)}, \quad \forall f\in L^q(\R),  
$$
which together with Sobolev's inequality and the fact that $\frac{1}{r}=\frac{1}{p}-\frac{\theta-1/2}{d}$ yield 
$$
\|\nabla u\|_{\mL^r_s} \leq C \|\p_t^{1-\theta} (-\Delta)^{\frac{1}{2}} u\|_{\mL^r_q}\leq C \|\p_t^{1-\theta} (-\Delta)^{\theta} u\|_{\mL^p_q}. 
$$
Combing this and \eqref{eq-frac}, we obtain 
\be\label{eq-pinter}
\|\nabla u\|_{\mL^r_s} \leq C\l( \|\p_t u\|_{\mL^p_q}+ \|\nabla^2 u\|_{\mL^{p}_q} \r),  
\ee
where $C$ only depends on $d, p, q, r, s$. 
If $u\in \mH^{2,p}_q(T)$, $\p_t u\in\mL^p_q(T)$ and $u(0,x)=0$, we extend $u$ as 
$$
\bar{u}(t,x):=\left\{\begin{array}{ll}
u(t, x) & \text { if } t\in [0, T]\\
-3 u\left(2T-t, x\right)+4 u\left(\frac{3T}{2}-\frac{t}{2}, x\right) & \text { if } t\in [T, 2T]\\
4 u\left(\frac{3T}{2}-\frac{t}{2}, x\right) & \text { if } t\in [2T, 3T]\\
0& \text { othewise. } 
\end{array}\right.
$$
By the definition of $\bar u$, one sees that 
$$
\|\nabla^k u\|_{\mL^p_q(T)}\asymp \|\nabla^k \bar u\|_{\mL^p_q}, \quad \|\p_t u\|_{\mL^p_q(T)}\asymp \|\p_t \bar u\|_{\mL^p_q}. 
$$
Therefore, our desired result follows from \eqref{eq-pinter}.  
\epf

\section{Kolmogorov's equation}
Throughout this paper, $Q$ always means a domain in $\R^{d+1}$ and $T>0$ is a time horizon. In this section, we study the unique solvability of the Kolmogorov equation \eqref{Eq-PDE} corresponding to \eqref{Eq-SDE} in some suitable $\widetilde \mH^{2, p_3}_{q_3}(T)$-space where $b$ satisfies the same assumptions as in Theorem \ref{Th-Main}. 

\subsection{Case 1: $b$ satisfies condition (a) of Theorem \ref{Th-Main}} 
\ 
\\
The main result in this subsection is 
\bt\label{Th-PDEmain1}
Let $d\geq 3$, $b=b_0+b_1$. 
Assume $b_1\in \mL^{p_1}_{q_1}(T)$ with  $\frac{d}{p_1}+\frac{2}{q_1}=1$ and $p_1\in (d, \infty)$. Then there for any  $p_3\in (1, d)$, $q_3\in (1, q_1)$ there is a constant $\eps=\eps(d, p_3, q_3)>0$ such that for each $f\in \widetilde \mL^{p_3}_{q_3}(T)$, equation \eqref{Eq-PDE} has a unique solution $u\in \widetilde \mH^{2,p_3}_{q_3}(T)$, provided that $\|b_0\|_{\tL^{d,\infty}_\infty(T)}\leq \eps$. Moreover, 
\be\label{Eq-W2est-1}
\|\p_t u\|_{\widetilde \mL^{p_3}_{q_3}(T)}+ \|u\|_{\widetilde \mH^{2,p_3}_{q_3}(T)} \leq C \|f\|_{\widetilde \mL^{p_3}_{q_3}(T)}, 
\ee
where $C$ only depends on $d, p_i, q_i, \eps, T$ and $b_1$. 
\et
\bpf

To prove the desired result, it suffices to show \eqref{Eq-W2est-1} assuming that the solution already exists, since the method of continuity is applicable. We first establish the corresponding estimate in the usual space $\mH^{2,p_3}_{q_3}(T)$, for any $p_3\in (1, d)$, $q_3\in (1,q_1)$ and some $\eps=\eps(d, p_3, q_3)>0$. 

Let $b^N_1:= b_1\1_{\{ |b_1|\leq N \}}$. Rewrite \eqref{Eq-PDE} as 
$$
\p_t u-\Delta u=f+ b_0\cdot\nabla u+ b_1^N\cdot\nabla u+(b_1-b_1^N)\cdot \nabla u. 
$$
Thanks to Theorem 1.2 of \cite{krylov2001heat}, for any $t\in [0,T]$, we have 
\const{\CW2}
\begin{align}\label{eq-I1}
\begin{aligned}
\|\p_t u\|&_{ \mL^{p_3}_{q_3}(t)} + \|\nabla^2 u\|_{ \mL^{p_3}_{q_3}(t)} \\
&\leq C_{\CW2}  \l( \|f\|_{\mL^{p_3}_{q_3}(t)}+ \|b_0\cdot \nabla u \|_{\mL^{p_3}_{q_3}(t)} + 
{N} \|\nabla u\|_{ \mL^{p_3}_{q_3}(t)}+
\|(b_1-b_1^N)\cdot \nabla u\|_{ \mL^{p_3}_{q_3}(t)} \r), 
\end{aligned}
\end{align}
where $C_{\CW2}=C_{\CW2}(d, p_3, q_3)$ does not depend on $t$. By \cite[Exercise 1.4.19]{grafakos2008classical} and  \cite[Remark 5]{tartar1998imbedding}, we get
\const{\CSobolev}  
\be\label{eq-I2}
\|b_0\cdot \nabla u \|_{\mL^{p_3}_{q_3}(t)} \leq C \|b_0\|_{\mL^{d,\infty}_{\infty}(t)}\|\nabla u\|_{\mL^{d{p_3}/(d-p_3), p_3}_{q_3}(t)}  \leq C_{\CSobolev} \|b_0\|_{\mL^{d,\infty}_{\infty}(t)}  \|\nabla^2 u\|_{\mL^{p_3}_{q_3}(t)}, 
\ee
where $\mL^{p, r}_q=L^q(\R_+; L^{q,r}(\R^d))$ ($L^{p,r}(\R^d)$ is the Lorentz space) and $C_{\CSobolev}=C_{\CSobolev}(d, p_3)$. Setting $1/r=1/p_3-1/p_1$ and $1/s=1/q_3-1/q_1$, by  \eqref{Eq-PSobolev1}, we have 
\begin{align}\label{eq-I3}
\begin{aligned}
\|(b_1-b_1^N)\cdot \nabla u\|_{ \mL^{p_3}_{q_3}(t)} \leq& \|(b_1-b_1^N)\|_{\mL^{p_1}_{q_1}(T)} \|\nabla u\|_{\mL^{r}_s(t)} \\
\leq& C_{\CPSobolevone} \|(b_1-b_1^N)\|_{\mL^{p_1}_{q_1}(T)} \l( \|\p_t u\|_{ \mL^{p_3}_{q_3}(t)} + \|\nabla^2 u\|_{ \mL^{p_3}_{q_3}(t)}\r). 
\end{aligned}
\end{align}
Let 
$$
\eps=\eps(d, p_3, q_3)=(4C_{\CW2}C_{\CSobolev})^{-1}>0.
$$ 
Noting that $q_1<\infty$, we can choose $N$ sufficiently large so that $\|b_1-b_1^N\|_{ \mL^{p_1}_{q_1}(T)}\leq (4C_{\CPSobolevone}C_{\CW2})^{-1}$.  By \eqref{eq-I1}-\eqref{eq-I3} and the choice of $\eps$ and $N$, if $\|b_0\|\leq \eps$, then for each $t\in [0,T]$, 
\const{\CIt}
\begin{align}\label{eq-It1}
I(t):= \|\p_t u\|_{ \mL^{p_3}_{q_3}(t)}^{q_3} +\|\nabla^2 u\|_{\mL^{p_3}_{q_3}(t)}^{q_3} \leq C_{\CIt}\l( \|f\|_{ \mL^{p_3}_{q_3}(t)}^{q_3}+ N^{q_3} \|\nabla u\|_{ \mL^{p_3}_{q_3}(t)}^{q_3} \r). 
\end{align}
Noting that 
\be\label{eq-u-I}
\begin{aligned}
\| u\|_{\mL^{p_3}_{q_3}(t)}^{q_3}=&\int_0^t \|u(s,\cdot)\|_{L^{p_3}}^{q_3} \d s=  
\int_0^t \l\|\int_0^s\p_r u(r,\cdot) \d r\r\|_{L^{p_3}}^{q_3} \d s\\
\leq & \int_0^t  s^{q_3-1} \|\p_t u\|_{ \mL^{p_3}_{q_3}(s)}^{q_3}  \d s \leq C(T, {q_3})\int_0^t I(s) \d s, 
\end{aligned}
\ee
using an interpolation inequality, we obtain 
\begin{align}\label{eq-inter}
\begin{aligned}
\|\nabla u\|_{ \mL^{p_3}_{q_3}(t)}^{q_3} \leq &\delta \|\nabla^2 u\|_{ \mL^{p_3}_{q_3}(t)}^{q_3} + C_\delta\| u \|_{ \mL^{p_3}_{q_3}(t)}^{q_3}\\
\leq &\delta I(t) + C_\delta \int_0^t I(s) \d s, \quad (\forall \delta>0). 
\end{aligned}
\end{align}
Combing \eqref{eq-It1} and \eqref{eq-inter}, we get 
\begin{align*}
I(t)  \leq & C_{\CIt} \delta N^{q_3} I(t)+C\|f\|_{ \mL^{p_3}_{q_3}(T)}^{q_3}+C_\delta N^{q_3} \int_0^t I(s) \d s.
\end{align*}
Letting $\delta=\delta(N)$ be small enough so that $ C_{\CIt} \delta N^{q_3}\leq 1/2$, we get 
$$
I(t)\leq C\|f\|_{ \mL^{p_3}_{q_3}(T)}^{q_3}+C \int_0^t I(s) \d s, \quad \forall t\in [0,T]. 
$$
Gronwall's inequality yields $I(T)\leq C \|f\|_{ \mL^{p_3}_{q_3}(T)}^{{q_3}}$, which together with \eqref{eq-u-I} implies 
\be\label{eq-W2-1}
\|\p_tu\|_{\mL^{p_3}_{q_3}(T)}+ \|u\|_{\mH^{2, p_3}_{q_3}(T)}\leq C \|f\|_{\mL^{p_3}_{q_3}(T)}, 
\ee
where $C=C(d, p_i, q_i, T, \eps, b_1)$. Our desired estimate \eqref{Eq-W2est-1} is then  obtained by \eqref{eq-W2-1} and an argument similar to the one in the proof for Lemma \ref {Le-heat}.
\epf

\subsection{Case 2: $b\in \tL^{d,\infty}_\infty(T)$ and $\div b\in \tL^{p_2}_\infty(T)$  with $p_2>d/2$.}  
\ 
\\
In this subsection, we will give an analogue of Theorem \ref{Th-PDEmain1}, where $b\in \widetilde \mL^{d,\infty}_\infty(T)$ and $\div b\in \tL^{p_2}_{\infty}(T)$. The result is stated as follows: 
\bt\label{Th-PDEmain2}
Let $d\geq 3$ and assume that $b\in \tL^{d,\infty}_\infty(T)$ and $\div b\in \tL^{p_2}_\infty(T)$ for some $p_2>d/2$.  Then there are constants  $p_3\in (d/2, d)$ and $q_3\in (2p_3/(2p_3-d), \infty)$ such that for each $f\in \tL^{d,\infty}_{\infty}(T)$, equation \eqref{Eq-PDE} has a unique solution $u\in \widetilde \mH^{2,p_3}_{q_3}(T)$. Moreover, 
\be\label{Eq-W2est-2}
\|\p_t u\|_{\widetilde \mL^{p_3}_{q_3}(T)}+ \|u\|_{\widetilde \mH^{2,p_3}_{q_3}(T)} \leq C \|f\|_{\widetilde \mL^{d,\infty}_{\infty}(T)}, 
\ee
where $C$ only depends on $d, p_2, p_3, q_3, T$ $\|b\|_{\tL^{d,\infty}_\infty(T)}$ and $\|\div b\|_{\tL^{p_2}_\infty(T)}$. 
\et

Unlike the previous case, if $\|b\|_{\widetilde \mL^{d,\infty}_{\infty}(T)}$ is large, then $\|b\cdot \nabla u\|_{\tL^{p_3}_{q_3}(T)}$ may not be controlled by\\
$\|\p_t u\|_{\tL^{p_3}_{q_3}(T)}+\|u\|_{\tH^{2, p_3}_{q_3}(T)}$, so the perturbation argument does not work any more. In order to overcome this difficulty, in this subsection,  by means of De Giorgi's method, we first show that any bounded weak solution of \eqref{Eq-PDE} is indeed  H\"older continuous, provided that $b$ is  in some Morrey's type space and $\div b\in \tL^{p_2}_\infty(T)$. Then in the light of Nirenberg's inequality \eqref{Eq-GNIeq2}, we show that $\nabla u$ is indeed in $\widetilde \mL^{r}_{q_3}(T)$ with some $r>d$, which implies $b\cdot\nabla u\in \widetilde \mL^{p_3}_{q_3}(T)$ with some $p_3>d/2$ and $q_3>2p_3/(2p_3-d)$. Our desired result then follows by Lemma \ref{Le-heat}. 

We first give the precise definition of weak solutions to the equation 
\be\label{Eq-PDE'}
\p_t u-\Delta u-b\cdot\nabla u=f \quad \mbox{in}\ Q=I\times D.  
\ee
\bd\label{Def-weak}
Assume $b\in L^2_{loc}(Q)$. We say $u\in V_{loc}(Q)$ is a subsolution (supersolution) to \eqref{Eq-PDE'} if 
for any $\varphi \in C_c^\infty(Q)$ with $\varphi\geq 0$,  
\be\label{Eq-weak}
\int_{Q} \big[-u \p_t \varphi +  \nabla u\cdot\nabla \varphi - b\cdot \nabla u \varphi \big]\leq(\geq ) \int_{Q} f \varphi.  
\ee
$u\in V_{loc}(Q)$ is a solution to \eqref{Eq-PDE'} if $u$ and $-u$ are subsolutions to \eqref{Eq-PDE'}. 
\ed

For any $p,q\in (1,\infty]$, here and below we define $p^*,q^*\in[2,\infty)$ by the relations
\begin{align}\label{Eq-*-pq}
\frac{1}{p}+\frac{2}{p^*}=1,\ \ \frac{1}{q}+\frac{2}{q^*}=1.
\end{align}
The following two lemmas are crucial for proving Theorem \ref{Th-PDEmain2}, and their proofs are essentially contained in \cite{zhang2020stochastic} and \cite{zhao2019stochastic}.  We provide sketches of their proofs in the Appendix for the reader's convenience.

\const{\Cenergy} 
\bl[Energy inequality]\label{Le-energy}
Assume $0<\rho<R\leq 1$, $k\geq 0$, $I\subseteq \mR$ is an open interval, $Q=I\times B_R$ and $\eta$ is a cut off function in $x$, compactly supported in $B_R$, $\eta(x)\equiv1$ in $B_{\rho}$, and $|\nabla\eta|\leq 2(R-\rho)^{-1}$.  Let  $d\geq 2$, $p_i, q_i\in (1,\infty)$ satisfying $d/p_i+2/q_i<2$, $i=2,3$. Suppose that $b, \div b \in\mL^{p_2}_{q_2}(Q)$, $f\in \mL^{p_3}_{q_3}(Q)$ and $u\in V(Q)$ is a bounded weak  subsolution to \eqref{Eq-PDE'}, then
\begin{align}\label{Eq-energy}
\begin{aligned}
&\l(\int u_k^2\eta^2\r)(t) - \l(\int u_k^2\eta^2 \r)(s) +\int_s^t\!\!\!\int |\nabla(u_k\eta)|^2\\
\leq &\frac{C_{\Cenergy}}{(R-\rho)^{2}}  \left(\|u_k\|^2_{\mL^2(A_s^t(k))}+ 
\sum_{i=2}^3\|u_k\|^2_{\mL^{p^*_i}_{q^*_i}(A_s^t(k))}\right)+C_{\Cenergy} \|f\|^2_{\mL^{p_3}_{q_3}(Q)} \|\1_{A_s^t(k)}\|_{\mL^{p^*_3}_{q^*_3}}^2, 
\end{aligned}
\tag{\bf{EI}}
\end{align}
where $u_k=(u-k)^+$, $A_s^t(k)= \{u>k\}\cap ([s,t]\times B_R)$ and $C_{\Cenergy}$ only depends on $d, p_i, q_i, \| b\|_{\mL_{q_2}^{p_2}(Q)}$  and $\|\div b\|_{\mL_{q_2}^{p_2}(Q)}$
\el

\const{\CLocalMax}
\bl\label{Le-Max}
Let $d\geq 2$ $p_i, q_i\in (1,\infty)$ satisfying $d/p_i+2/q_i<2$, $i=2,3$. Suppose $b, \div b \in\mL^{p_2}_{q_2}(Q_1)$ and
$u\in V(Q_1)$ is a bounded weak subsolution to \eqref{Eq-PDE'} in $Q_1$. Then for any $f\in \mL^{p_3}_{q_3}(Q_1)$, 
\begin{align}\label{Eq-local-max}
\begin{aligned}
\|u^+\|_{L^\infty(Q_{1/2})} \leq &C_{\CLocalMax} \left(\|u^+\|_{\mL^2_2(Q_1)}+\sum_{i=2}^3\|u^+\|_{\mL^{p^*_i}_{q^*_i}(Q_1)}+ \|f\|_{\mL_{q_3}^{p_3}(Q_1)}\right). 
\end{aligned}
\tag{ \bf{LM} }
\end{align}
Here $C_{\CLocalMax}$ only depends on $d, p_i, q_i, \| b\|_{\mL_{q_2}^{p_2}(Q_1)}$  and $\|\div b\|_{\mL_{q_2}^{p_2}(Q_1)}$.
\el

\medskip

In order to prove the H\"older estimate for the bounded weak solutions to \eqref{Eq-PDE}, we also need some technical  Lemmas. One of them is a parabolic version of De Giorgi's Lemma.  Set $Q'_1:= (-2,-1)\times B_1.$ For any $u: Q_1\cup Q_1'\to \R$, define  
$$
A_u:=\{f\geq 1/2\}\cap Q_1, \quad  B_u:= \{f\leq 0\}\cap Q'_1, \quad D_u:=\{0<f<1/2\}\cap (Q_1\cup Q'_1). 
$$
\bl\label{Le-parabolic-DeGiorgi}
Let $p_i, q_i\in (1,\infty)$ satisfying $d/p_i+2/q_i<2$, $i=2,3$. Assume $b, \div b \in \mL^{p_2}_{q_2}(Q_2)$, $\|f\|_{\mL^{p_3}_{q_3}(Q_2)}\leq 1$ and that $u$ is a weak subsolution  to \eqref{Eq-PDE'}  in $Q_2$ with  $u\leq 1$.  
Suppose that $\delta\in (0,1)$ and that 
\begin{align*}
|A_{u}|\geq \delta\ \mbox{ and }\ |B_{u}|\geq \delta. 
\end{align*}
Then 
$$
|D_u|=|\{0<u<1/2\}\cap (Q_1\cup Q'_1)|\geq \beta, 
$$ 
where $\beta=\beta(d, p_i,q_i,\|b\|_{\mL^{p_2}_{q_2}(Q_2)}, \|\div b\|_{\mL^{p_2}_{q_2}(Q_2)}, \delta)$ is a universal constant that does not depend on $u$.  
\el
\begin{proof}
By \eqref{Eq-energy}, H\"older's inequality and our assumption $u\leq 1$,  we have 
\begin{align}\label{Eq-energy-u+}
\begin{aligned}
&\l( \int_{B_1} (u^+)^2 (t)-\int_{B_1} (u^+)^2 (s)\r)+ \int_s^t\!\!\!\int_{B_1} |\nabla (u^+)|^2 \d x\d r\\
\leq&  C \left(\|u^+\|_{\mL^{2}_{2}([s,t]\times B_2)}^2+ \sum_{i=2}^3\|u^+\|_{\mL^{p^*_i}_{q^*_i}([s,t]\times B_2)}^2 + \|f\|^2_{\mL_{q_3}^{p_3}(Q_2)} \|\1_{\{[s,t]\times B_2\}}\|_{\mL^{ p^*_3}_{q^*_3}}^2 \right) \\
\leq & C |t-s|^{\theta},
\end{aligned}
\end{align}
where $\theta=\frac{1}{2}\wedge (1-\frac{1}{q_2}) \wedge (1-\frac{1}{q_3})>0$. Assume $|D_u|<\beta$, where $\beta>0$ is a small number, which will be determined later. Let 
\begin{align*}
&a(t)= |\{x\in B_1: u^+(t, x) \geq 1/2\}|, \\
&b(t)= |\{x\in B_1: u^+(t, x) = 0\}|, \\
&d(t)= |\{x\in B_1: 0<u^+(t, x)<1/2\}|.
\end{align*} 
Set
$$
I_1:= \{t\in (-2,0): d(t)\leq \sqrt{\beta}\}\ \mbox{ and }\ I_2:= \l\{t\in I_1 :   b(t)>|B_1|-\frac{\delta}{100d!} \mbox{ or } b(t)<\frac{\delta}{100d!}  \r\}. 
$$  
By our assumption and Chebyshev's inequality, $|I_1|\geq 2-C\sqrt{\beta}$. Using \eqref{Eq-energy-u+} and Lemma \ref{Le-De Giorgi}, we have 
$$
C\geq \int_{I_1}\!\int_{B_1}|\nabla u^+(t,x)|^2 \d x \geq c_d\beta^{-\frac{1}{2}}\int_{I_1} a^2(t) b^{2-\frac{2}{d}}(t) \d t. 
$$
Thus, 
$$
\int_{I_1} a^2(t) b^{2-\frac{2}{d}} (t) \d t \leq C\sqrt \beta \to 0 \mbox{ as } \beta\to 0. 
$$
This together with the facts that $\inf_{t\in I_1} [a(t)+b(t)]\geq |B_1|-\sqrt{\beta}\geq 1/d!$ and $|I_1|\to 2$ implies $|I_2|\to 2$ as $\beta\to 0$. Since the zero set of $u^+$ has mass $\delta$ in $Q_1'$, for small $\beta$, there is some $t_1\in (-2,-1)\cap I_2$ such that $b(t_1)>|B_1|-\frac{\delta}{100d!}$. Using the first term in the energy estimate \eqref{Eq-energy-u+}, 
we see that for some universal small $\tau>0$, there exists $t_2\in I_2$ with $t_2\geq t_1+\tau$ such that for all $t\in [t_1, t_2]\cap I_2$, $b(t)>|B_1|-\frac{\delta}{100d!}$. Iterating, we obtain that for all $t\in [t_1, 0]\cap I_2$, $a(t)<\frac{\delta}{100d!}$. That is a contradiction to $|A_u|\geq \delta$ and completes the proof.
\end{proof}

The next diminish of oscillation lemma is crucial for the H\"older estimate for the solutions to \eqref{Eq-PDE}.  
\bl\label{Le-diminish}
Let $p_i, q_i\in (1,\infty)$ satisfying $d/p_i+2/q_i<2$ , $i=2,3$. Assume $b\in \mL^{p_2}_{q_2}(Q_2)$ and $\div b\in \mL^{p_2}_{q_2}(Q_2)$. Then there exist universal positive constants $\mu<1$ and $\eps_0>0$ only depending on $d, p_i,q_i,\|b\|_{\mL^{p_2}_{q_2}(Q_2)}$ and $\|\div b\|_{\mL^{p_2}_{q_2}(Q_2)}$, such that for any $f\in\mL^{p_3}_{q_3}(Q_2)$ with $\|f\|_{\mL^{p_3}_{q_3}(Q_2)}\leq \eps_0$, and any weak subsolution $u$ of \eqref{Eq-PDE'} in $Q_2$ satisfying $u \leq 1$ and $|\{u\leq 0\}\cap Q'_1|\geq |Q'_1|/2$, the following estimate is valid: 
$$
u\leq \mu \ \mbox{ in } Q_{1/2}. 
$$
\el
\begin{proof}
We consider $u_k = 2^k(u-(1-2^{-k}))$, which fulfills for each $k\geq 0$, 
$$u_k\leq 1, \quad B_{u_k}:=|\{u_k\leq 0\}\cap Q'_1| \geq |Q'_1|/2
$$ 
and 
$$
\p_t u_k-L u_k=2^k f=:f_k. 
$$  
Let $\delta\in (0,1)$ be sufficiently small such that  $8 C_{\CLocalMax}\delta\leq 1$, where $C_{\CLocalMax}$ is the constant in \eqref{Eq-local-max}. Suppose $\beta=\beta_{\delta}>0$ is the same constant as in Lemma \ref{Le-parabolic-DeGiorgi} and set 
$$
K:=[3|B_1|/\beta]+1\ \mbox{ and }\ \eps_0:= 2^{-K} \delta. 
$$ 
By the definitions of $K$ and $\eps_0$, one can see that for each $k\in \{1,2,\cdots, K\}$,  
\be\label{eq-fk}
\|f_k\|_{\mL^{p_3}_{q_3}(Q_2)}=2^k\| f\|_{\mL^{p_3}_{q_3}(Q_2)}\leq 2^K \eps_0= \delta<1.
\ee
We claim that 
\be\label{eq-uKless}
\|u^+_{K}\|^2_{\mL^2_2(Q_1)}+\sum_{i=2}^3
\|u^+_{K}\|^2_{\mL^{p^*_i}_{q^*_i}(Q_1)} \leq 3\delta. \ee
Assume \eqref{eq-uKless} does not hold. Noting that $u_k$ is decreasing, so  
\be\label{eq-uklarger}
\|u_{k}^+\|^2_{\mL^2_2(Q_1)}+\sum_{i=2}^3
\|u^+_{k}\|^2_{\mL^{p^*_i}_{q^*_i}(Q_1)}>3 \delta, 
\ee
for all $k\in \{1,2,\cdots, K\}$. By \eqref{eq-uklarger} and the fact that $u_k\leq 1$, we get
$$
\l|\{u_{k}\geq 0\} \cap Q_1\r| \geq \frac13\l(\|u_{k}^+\|^2_{\mL^2_2(Q_1)}+\sum_{i=2}^3
\|u^+_{k}\|^2_{\mL^{p^*_i}_{q^*_i}(Q_1)} \r)>\delta. 
$$
Thus, 
\begin{align*}
|A_{u_{k-1}}|:=&\{u_{k-1}\geq 1/2\}\cap Q_1|= |\{u_{k}\geq 0\}\cap Q_1|>\delta.  
\end{align*}
Recalling that $u_{k-1}\leq 1$, $|B_{u_{k-1}}|\geq |Q'_1|/2$ and $\|f_{k-1}\|_{\mL^{p_3}_{q_3}(Q_2)}\leq 1$, by virtue of Lemma \ref{Le-parabolic-DeGiorgi}, we have 
$$
|\{1-2^{-k+1}<u<1-2^{-k}\}\cap (Q_1\cup Q'_1)|=|D_{u_{k-1}}| \geq \beta.  
$$
Hence, 
\begin{align*}
2|B_1|\geq& |\{0<u<1-2^{-K}\} \cap (Q_1\cup Q'_1)|\\
\geq &\sum_{k=1}^K |\{1-2^{-k+1}<u<1-2^{-k}\}\cap (Q_1\cup Q'_1)|\\
\geq &K\beta=([3|B_1|/\beta]+1)\beta\geq 3|B_1|,  
\end{align*}
which is a contradiction. So we complete the proof for \eqref{eq-uKless}. This together with \eqref{Eq-local-max} yields 
\begin{align*}
\|u^+_{K}\|_{L^\infty(Q_{1/2})}\leq & C_{\CLocalMax}\left(\|u^+_{K}\|^2_{\mL^2_2(Q_1)}+\sum_{i=2}^3
\|u^+_{K}\|^2_{\mL^{p^*_i}_{q^*_i}(Q_1)}+\|f_K\|_{\mL^{p_3}_{q_3}(Q_1)}\right)\\
\leq &  4C_{\CLocalMax}\delta=1/2,  
\end{align*}
which implies 
$$
\sup_{x\in Q_{1/2}} u \leq 1-2^{-K-1}. 
$$
Letting $\mu=1-2^{-K-1}$, we complete our proof.  
\end{proof}
\bl
Let $p_i, q_i\in (1,\infty)$ satisfy $d/p_i+2/q_i<2$ , $i=2,3$. Assume $b, \div b\in \mL^{p_2}_{q_2}(Q_2)$. Suppose $u$ is a weak solution to \eqref{Eq-PDE'} in $Q_2$ with $f\in \mL^{p_3}_{q_3}(Q_2)$. Then 
\be\label{eq-osc}
\osc_{Q_{1/2}} u\leq \mu\, \osc_{Q_{1}} u+C  \|f\|_{\mL_{q_3}^{p_3}(Q_2)}, 
\ee
where $\mu<1$ is the same constant as in Lemma \ref{Le-diminish} and $C$ only depends on $d, p_i,q_i$,$\|b\|_{\mL^{p_2}_{q_2}{(Q_2)}}$ and $\|\div b\|_{\mL^{p_2}_{q_2}(Q_2)}$. 
\el
\bpf
Define 
$$
u_\delta:=\frac{u}{\delta+\|u^+\|_{L^\infty(Q_2)}+\eps_0^{-1}\|f\|_{\mL^{p_3}_{q_3}(Q_2)}}\leq 1,\quad \delta>0, 
$$
where $\eps_0$ is the same constant as in Lemma \ref{Le-diminish}. Then $u_\delta$ satisfies 
$$
\p_t  u_\delta-L u_\delta=f_\delta:= f(\delta+\|u^+\|_{L^\infty(Q_2)}+\eps_0^{-1}\|f\|_{\mL^{p_3}_{q_3}(Q_2)})^{-1}
$$
in $Q_2$ and $\| f_\delta\|_{\mL^{p_3}_{q_3}(Q_2)}\leq \eps_0$. By Lemma \ref{Le-diminish}, we have $u_\delta\leq \mu$ in $Q_{1/2}$, so 
\begin{align*}
\|u^+\|_{L^\infty(Q_{1/2})} \leq & \mu \liminf_{\delta\downarrow 0} \left(\delta+\|u^+\|_{L^\infty(Q_2)}+\eps_0^{-1}\|f\|_{\mL^{p_3}_{q_3}(Q_2)}\right)\\
\leq & \mu \|u^+\|_{L^\infty(Q_{2})} + \mu \eps_0^{-1} \|f\|_{\mL^{p_3}_{q_3}(Q_2)}. 
\end{align*}
Similarly, $\|u^-\|_{L^\infty(Q_{1/2})} \leq \mu \|u^-\|_{L^\infty(Q_{2})} + \mu \eps_0^{-1} \|f\|_{\mL^{p_3}_{q_3}(Q_2)}.$ So, we complete our proof. 
\epf
Now we are at the point to show the H\"older regularity of the solutions to \eqref{Eq-PDE}. 

\bl\label{Le-Holder}
Let $p_i, q_i\in (1,\infty)$ satisfying $d/p_i+2/q_i<2$ , $i=2,3$. Suppose that $b, \div b \in\widetilde \mL^{p_2}_{q_2}(T)$   and $f\in \widetilde\mL^{p_3}_{q_3}(T)$. If there is a constant $\cK_b<\infty$ such that for all $(t,x) \in [0,T]\times\R^d$ and $r\in (0,1)$,  
$$
r^{1-\frac{d}{p_2}-\frac{2}{q_2}}  \|b\|_{\mL^{p_2}_{q_2}(Q_{2r}(t,x))}\leq \cK_b,
$$ 
then there are constants $\a\in (0,1)$ and $C>1$ such that for any bounded weak solution $u\in \widetilde V^{0}(T)$ to \eqref{Eq-PDE}, it holds that 
\be\label{Eq-Holder}
\|u\|_{C^\a([0,T]\times\R^d)} \leq C \|f\|_{\widetilde \mL^{p_3}_{q_3}(T)}, 
\ee
where $\a, C$ only depend on $d, p_i, q_i, T, \|b\|_{\widetilde \mL^{p_2}_{q_2}(T)}$ and $\cK_b$. 
\el

\begin{proof}
For convenience, we extend $u, b, f$ to be functions on $(-\infty, T)\times \R^d$ by letting $u(t,x)=b(t,x)=f(t,x)=0$, for all $t\leq 0$ and $x\in \R^d$. By Definition \ref{Def-weak}, $u$ is still a bounded weak solution to \eqref{Eq-PDE'} on $(-\infty, T)\times \R^d$. For any $r\in (0,1)$, $(t_0, x_0)\in (-\infty, T)\times \R^d$ and $(t,x)\in Q_2$, define $u_r(t,x):= u(r^2 t+t_0, rx+x_0)$, $b_r(t,x):= rb(r^2 t+t_0, rx+x_0)$, $f_r(t,x):= r^2f(r^2 t+t_0, rx+x_0)$. Then $u_r$ satisfies 
$$
\p_t u_r -\Delta u_r -b_r\cdot\nabla u_r=f_r \ \mbox{ in } \ Q_2. 
$$
By our assumption, we have  
\begin{align*}
&\|b_r\|_{\mL^{p_2}_{q_2}(Q_{2})} = r^{1-\frac{d}{p_2}-\frac{2}{q_2}} \|b\|_{\mL^{p_2}_{q_2}(Q_{2r}(t_0, x_0))}\leq \cK_b,\\
&\|\div b_r\|_{\mL^{p_2}_{q_2}(Q_{2})} =r^{\kappa_2} \|\div b\|_{\mL^{p_2}_{q_2}(Q_{2r}(t_0, x_0))}\leq \|\div b\|_{\tL^{p_2}_{q_2}(T)}, \\
&\|f_r\|_{\mL^{p_3}_{q_3}(Q_2)}= r^{\kappa_3} \|f\|_{\mL_{q_3}^{p_3}(Q_{2r}(t_0, x_0))}, 
\end{align*}
where $\kappa_i= 2-d/p_i-2/q_i >0$, $i=2,3$. Using \eqref{eq-osc}, we get 
\be\label{eq-osc2}
\osc_{Q_{\frac{r}{2}}(t_0,x_0)} u \leq \mu \osc_{Q_{2r}(t_0,x_0)} u + C r^{\kappa_3} \|f\|_{\mL_{q_3}^{p_3}(Q_{2r}(t_0, x_0))}, \ \mu\in (0,1).   
\ee
The desired estimate \eqref{Eq-Holder} follows by \eqref{eq-osc2} and standard arguments (see \cite[Lemma 3.4]{han2011elliptic}). 
\end{proof}

\br
We should also point out that the Harnack inequality for Lipschitz continuous solutions to \eqref{Eq-PDE'} with $f\equiv0$ was also obtained in \cite{nazarov2012harnack}  by Moser iteration method. 
\er

To prove our desired result, we also need the following simple lemma. 
\bl\label{Le-Lp-Lrw}
Let $1<p<r<\infty$ and $A$ be a Borel subset of $\R^d$ with finite Lebesgue measure. Then, there is a constant $C=C(d, p, r)$ such that 
\be\label{Eq-Lp-wLd}
\|f\|_{\widetilde L^p} \leq C(d, p, r) \|f\|_{\widetilde L^{r,\infty}}. 
\ee
\el
\bpf
Let $A$ be any Borel subset of $\R^d$. Set 
$$
\mu_f(t)= |\{x\in A: |f(x)|>t\}|. 
$$
Then, 
\begin{align*}
\int_{A} |f|^p =&p\int_0^\infty t^{p-1} \mu_f(t) \d t=p\int_0^{\lambda} t^{p-1} |A| \d t + p\|f\|_{L^{r,\infty}(A)}^r \int_{\lambda}^\infty t^{p-r-1} \d t\\
\leq &  \lambda^p |A| +p (r-p)^{-1} \|f\|_{L^{r,\infty}(A)}^r \lambda^{p-r}. 
\end{align*}
Letting $\lambda=(\frac{p}{r-p})^{1/r}\|f\|_{L^{r,\infty}(A)} |A|^{-1/r}$, we obtain 
$$
\|f\|_{L^p(A)} \leq 2^{1/p} \left(\frac{p}{r-p}\right)^{1/r} \|f\|_{L^{r,\infty}(A)} |A|^{1/p-1/r}. 
$$
Thus, 
$$
\|f\|_{\widetilde L^p} \leq \sup_{y\in \R^d} \|f\|_{L^p(B_2(y))} \leq C(d, p, r) \|f\|_{\widetilde L^{r,\infty}}. 
$$
\epf

Now we are in the position of proving Theorem  \ref{Th-PDEmain2}. 
\bpf[Proof of Theorem \ref{Th-PDEmain2}] 
Since $\widetilde L^{p}\subseteq \widetilde L^{p'}\, (p>p')$,  we can assume $p_2\in (d/2, d)$. Letting $q_2\in (1,\infty)$ such that $d/p_2+2/q_2<2$, by our assumptions on $b$, one sees that $b, \div b\in \widetilde \mL^{p_2}_{q_2}(T)$ and  for any $r\in (0,1)$, $t_0\in [0,T]$ and $x_0\in \R^d$, 
\begin{align*}
&r^{1-\frac{d}{p_2}-\frac{2}{q_2}}\l(\int_{t_0-r^2}^{t_0}\|b(t,\cdot)\|_{L^{p_2}(B_r(x_0))}^{q_2} \ \d t\r)^{1/q_2} \\
\overset{\eqref{Eq-Lp-wLd}}{\leq} &C r^{-\frac{2}{q_2}}\l( \int_{t_0-r^2}^{t_0}\|b(t, \cdot)\|_{L^{d, \infty}(B_r(x_0))}^{q_2}\ \d t \r)^{1/q_2}\leq C  \|b\|_{\widetilde \mL^{d, \infty}_\infty(T)}. 
\end{align*}
Let $p_3'\in (d/2, d), q_3'\in (1,\infty)$ be some constants satisfying $d/p_3'+2/q_3'<2$. Again by Lemma \ref{Le-Lp-Lrw}, $f\in \tL^{p_3'}_{q_3'}(T)$. Thanks to Lemma \ref{Le-GlobalMax} and Theorem \ref{Le-Holder}, \eqref{Eq-PDE} admits a unique bounded weak solution $u$, and there is a constant $\a\in (0, 1)$ only depending on $d, p_2, q_2, p_3', q_3', T, \|b\|_{\tL^{p_2}_{q_2}(T)}$ and $\|\div b\|_{\tL^{p_2}_{q_2}(T)}$ such that  
\be\label{eq-holder}
\|u\|_{C^\a([0,T]\times \R^d)} \leq C \|f\|_{\widetilde \mL^{p_3'}_{q_3'}(T)} \leq C \|f\|_{\tL^{d, \infty}_\infty(T)}. 
\ee
Next we fix 
\be\label{eq-s-q1}
s\in \l(2\vee \frac{d(4-3\a)}{4-2\a}, d\r), \quad q_3>\frac{2s(2-\a)}{2s(2-\a)-d(4-3\a)}.  
\ee
Rewrite \eqref{Eq-PDE} as 
$$
\p_t u-\Delta u= f+\div (b u)-(\div b) u. 
$$
It is easy to see that 
$$
\|f\|_{\widetilde \mH^{-1,s}_{q_3}(T)} \leq C \|f\|_{\widetilde \mL^{d, \infty}_{\infty}(T)}, \quad \|\div (b u)\|_{\widetilde \mH^{-1,s}_{q_3}(T)} \leq C \|b\|_{\widetilde \mL^s_{q_3}(T)} \|u\|_{\mL^\infty(T)}\overset{\eqref{eq-holder}}{\leq} C\|f\|_{\widetilde \mL^{d, \infty}_{\infty}(T)}. 
$$
Noting that $sd/(d+s)<d/2<p_2$ and using Sobolev embedding, we see that
\begin{align*}
\|(\div b) u\|_{\widetilde \mH^{-1,s}_{q_3}(T)} \leq& C \|(\div b) u\|_{\tL^{sd/(s+d)}_{q_3}(T)}\\
\leq& C \|\div b\|_{\tL^{p_2}_{\infty}(T)} \|u\|_{\mL^\infty(T)}\overset{\eqref{eq-holder}}{\leq} C\|f\|_{\widetilde \mL^{d, \infty}_{\infty}(T)}. 
\end{align*}
By Lemma \ref{Le-heat}, 
$$
\|u\|_{\widetilde \mH^{1,s}_{q_3}(T)} \leq C \|f+\div (b u)-(\div b) u\|_{\widetilde \mH^{-1,s}_{q_3}(T)}\leq C \|f\|_{\widetilde \mL^{d, \infty}_{\infty}(T)}. 
$$
Using this and noting the fact that $1<s/2<d/2$, we get  
$$
\|b\cdot \nabla u\|_{\widetilde \mL^{s/2}_{q_3}(T)} \leq \|b\|_{\widetilde \mL^{d, \infty}_\infty(T)} \|u\|_{\widetilde \mH^{1,s}_{q_3}(T)} \leq C\|f\|_{\widetilde \mL^{d, \infty}_{\infty}(T)}. 
$$
Again by Lemma \ref{Le-heat}, we obtain  
\be\label{eq-u''}
\|u\|_{\widetilde \mH^{2, s/2}_{q_3}(T)} \leq C\|f+b\cdot\nabla u\|_{\widetilde \mL^{s/2}_{q_3}(T)}\leq C\|f\|_{\widetilde \mL^{d, \infty}_{\infty}(T)}. 
\ee
In the light of  Nirenberg's inequality \eqref{Eq-GNIeq2}, we get
$$
\|\nabla u\|_{\widetilde \mL^{r}_{q_3}(T)} \overset{\eqref{Eq-GNIeq2}}{\leq} C \|\nabla^2 u\|_{\widetilde \mL^{s/2}_{q_3}(T)}^\theta \cdot \|u\|_{C^\a([0,T]\times \R^d)}^{1-\theta}\overset{\eqref{eq-holder}, \eqref{eq-u''}}{\leq} C\|f\|_{\widetilde \mL^{d, \infty}_{\infty}(T)}, 
$$ 
where 
\be\label{eq-r}
r=\frac{(2-\a)s}{2-2\a}, \quad \theta= \frac{s}{2r}\in (0,1). 
\ee 
Now letting  
\be\label{eq-p1}
\frac{1}{p_3}=\frac{1}{r}+\frac{1}{s}= \frac{4-3\a}{s(2-\a)}, 
\ee
by H\"older's inequality, 
$$
\| b\cdot \nabla u\|_{\widetilde \mL^{p_3}_{q_3}(T)}\leq C \|b\|_{\widetilde \mL^s_\infty(T)} \|\nabla u\|_{\widetilde \mL^{r}_{q_3}(T)} \leq C \|b\|_{\widetilde \mL^{d, \infty}_\infty(T)} \|f\|_{\widetilde \mL^{d, \infty}_{\infty}(T)}. 
$$
Using Lemma \ref{Le-heat} again, we obtain   
$$
\|\p_t u\|_{\widetilde \mL^{p_3}_{q_3}(T)}+\|u\|_{\widetilde \mH^{2,p_3}_{q_3}(T)} \leq C \|f\|_{\widetilde \mL^{d, \infty}_{\infty}(T)}. 
$$
By \eqref{eq-s-q1}, \eqref{eq-r} and \eqref{eq-p1},  we have 
$$
\frac{d}{p_3}+\frac{2}{q_3}= \frac{d(4-3\a)}{s(2-\a)}+\frac{2}{q_3}<2. 
$$
So we complete our poof. 
\epf

\section{Proof of the main result}
In this section, we present the proof of our main probabilistic result. Firstly, we give the precise definition of martingale solutions to \eqref{Eq-SDE}. 
\bd
For given $x\in\mR^d$, we call a probability measure $\mP_{x}\in\sP(C([0,T]; \R^d))$ a martingale solution of SDE \eqref{Eq-SDE} with starting point $x$ if
\begin{enumerate}[(i)]
\item $\mP_{x}(\omega_0=x)=1$, and for each $t\in [0,T]$,
$$
\mE_x \int^t_0 |b(s,\omega_s)|\dif s<\infty,  
$$
where $\{\om_t\}_{t\in [0,T]}$ is the canonical processes.  
\item For all $f\in C^2_c(\mR^d)$, 
$$
M^f_t(\omega):=f(\omega_t)-f(x)-\int^t_0 \l(\Delta f+b\cdot\nabla f\r) (\omega_s)\dif s
$$
is a $\cB_t$-martingale under $\mP_{x}$, where $\cB_t:= \sigma\left\{\omega_{s}: 0 \leqslant s \leqslant t\right\}$. 
\end{enumerate}
\ed

Let $\rho\in C_c^\infty(\R^d)$ and $\int_{\R^d} \rho=1$. Set $\rho_n(x):= n^d \rho(nx)$ and $b_n(t,x)=b(t,\cdot)*\rho_n(x)$. For each $x\in\mR^d$, we then consider the following modified SDE:
\begin{align}\label{SDE9}
\dif X^n_{t}(x)=b_n(t,X^n_{t}(x))\dif t+\sqrt{2}\dif W_t,\ \ X^n_{0}=x,
\end{align}
where $W$ is a $d$-dimensional standard Brownian motion on some complete filtered probability space 
$(\Omega,\sF,(\sF_t)_{t\in [0,T]},\bP)$. 
It is well known that there is a unique strong solution $X^n_{t}(x)$ to the above SDE.

\bpf[Proof of Theorem \ref{Th-Main}]
{\em Existence:}
Assume $b$ satisfies condition (a) or (b) in Theorem \ref{Th-Main} and $p_3, q_3\in (1,\infty)$ such that $d/p_3+2/q_3<2$. We first prove that there are constants $\theta>0$ 
and $C>0$ such that for any $f\in C^\infty_c(\mR^{d+1})$
and $0\leq  t_0<t_1\leq  T$, 
\begin{align*}
\sup_n\sup_{x\in\mR^d}\bE \int^{t_1}_{t_0} f(t,X_{t}^n(x))\dif t \leq C(t_1-t_0)^\theta\| f\1_{[t_0,t_1]}\|_{\widetilde \mL^{p_3}_{q_3}}.
\end{align*}
Let $u_n$ be the smooth solution of the following backward PDE:
\begin{align}\label{eq-un}
\p_t u_n+\Delta u_n+b_n\cdot\nabla u_n+f=0,\ u_n(t_1,\cdot)=0.
\end{align}
By It\^o's formula we have
$$
u_n(t_1,X^n_{t_1})=u_n(t_0,X^n_{t_0})+\int^{t_1}_{t_0}
(\p_t u_n+\Delta u_n+b_n\cdot\nabla u_n)(t,X^n_{t})\dif t+\sqrt{2}\int^{t_1}_{t_0}\nabla u_n(t,X^n_{t})\dif W_t.
$$
Using \eqref{eq-un} and taking expectation, we obtain
\begin{align*}
\bE \int^{t_1}_{t_0}f(t,X^n_{t})\dif t =\bE u_n(t_0,X^n_{t_0})| \leq \|u_n(t_0,\cdot)\|_{L^\infty}.
\end{align*}
Since $\frac{d}{p_3}+\frac{2}{q_3}<2$, we can choose $q'_3<q_3$ so that $\frac{d}{p_3}+\frac{2}{q'_3}<2$.
Thus, by Lemma \ref{Le-GlobalMax}, we obtain 
\be\label{eq-kry1}
\bE \int^{t_1}_{t_0}f(t,X^n_{t})\dif t\leq \|u_n(t_0,\cdot)\|_{L^\infty} 
\leq C\| f\1_{[t_0,t_1]}\|_{\widetilde \mL^{p_3}_{q'_3}}\leq C(t_1-t_0)^{1-\frac{q'_3}{q_3}}\|f\1_{[t_0,t_1]}\|_{\widetilde \mL^{p_3}_{q_3}}.
\ee
Now let $\tau\leq T$ be any bounded stopping time. Note that
$$
X^n_{(\tau+\delta)\wedge T}(x)-X^n_{\tau}(x)=\int^{(\tau+\delta)\wedge T}_\tau b_n(t,X^n_{t}(x))\dif t+\sqrt{2}(W_{(\tau+\delta)\wedge T}-W_\tau),\ \ \delta\in (0, 1).
$$
By \eqref{eq-kry1} and Remark 1.2 in \cite{zhang2020stochastic}, we have
$$
\bE \int^{(\tau+\delta)\wedge T}_\tau |b_n|(t,X^n_{t}(x))\dif t \leq C\delta^\theta\| b_n\|_{\widetilde \mL^{p_3}_{q_3}(T)}. 
$$
Thus, 
\begin{align*}
\bE\sup_{0\leq u\leq \delta}|X^n_{\tau+u}(x)-X^n_{\tau}(x)|&
\leq \bE \int^{\tau+\delta}_\tau |b_n|(t,X^n_{t}(x))\dif t+\sqrt{2} \bE \sup_{0\leq u\leq \delta}|W_{\tau+\delta}-W_\tau|\\
&{\leq } C\delta^\theta\| b_n\|_{\widetilde \mL^{p_3}_{q_3}}+C\delta^{1/2}\leq C\delta^{\theta'},
\end{align*}
where $\theta'>0$ and $C$ is independent of $n$.  So by \cite[Lemma 2.7]{zhang2018singular}, we obtain
$$
\sup_n\sup_{x\in \mR^d}\bE\left(\sup_{t\in[0,T]; u\in [0,\delta]}|X^n_{t+u}(x)-X^n_{t}(x)|^{1/2}\right)\leq C\delta^{\theta'}.
$$
From this, by Chebyshev's inequality, we derive that for any $\eps>0$,
$$
\lim_{\delta\to 0}\sup_n\sup_{x\in\mR^d}\bP\left(\sup_{t\in[0,T]; u\in [0,\delta]}|X^n_{t+u}(x)-X^n_{t}(x)|>\eps\right)=0.
$$
Hence, by \cite[Theorem 1.3.2]{stroock2007multidimensional}, $\mP^n_x:=\bP \circ X^n_{\cdot}(x)^{-1}$ is tight in $\sP(C([0,T]; \R^d))$. Assume $\mP_{x}$ is an accumulation point of $(\mP^n_{x})_{n\in\mN}$, that is, for some subsequence $n_k$,
$$
\mP^{n_k}_{x}\mbox{ weakly converges to $\mP_{x}$ as $k\to\infty$. }
$$
Since \eqref{eq-kry1} can be rewritten as 
$$
\mE^n_{x}\left(\int^{t_1}_{t_0} f(t,\omega_t)\dif t\right)\leq C(t_1-t_0)^\theta\| f \1_{[t_0,t_1]}\|_{\widetilde \mL^{p_3}_{q_3}}. 
$$
By taking weak limits and a standard monotone class argument, we obtain 
\begin{align}\label{eq-kry2}
\mE_{x}\left(\int^{t_1}_{t_0} f(t,\omega_t)\dif t\right)\leq C(t_1-t_0)^\theta\| f \1_{[t_0,t_1]}\|_{\widetilde \mL^{p_3}_{q_3}}, \quad \mbox{ for all } d/p_3+2/q_3<2. 
\end{align}
In order to prove that $\mP_x$ is a martingale solution to \eqref{Eq-SDE}, it suffices to prove that for any $0\leq t_0<t_1\leq T$ and $f\in C^2_c(\mR^d)$,
$$
\mE_{x}(M^f_{t_1}|\cB_{t_0})=M^f_{t_0},\ \ \ \mP_{x}-a.s.,
$$
where 
$$
M^f_t:=f(\omega_t)-
f(\om_0)-\int^t_0(\Delta+b\cdot\nabla) f(s,\omega_s)\dif s.
$$
By a standard monotone class argument, it is enough to show that for any $G\in C_b(C([0,T]; \R^d))$ being $\cB_{t_0}$-measurable, 
$$
\mE_{x}\Big(M^f_{t_1}\cdot G\Big)=\mE_{x}\Big(M^f_{t_0}\cdot G\Big).
$$
Note that for each $n\in\mN$,
$$
\mE^n_{x}\Big(
M^{n,f}_{t_1}\cdot G\Big)
=\mE^n_{x}\Big(
M^{n,f}_{t_0}\cdot G\Big), 
$$
where 
$$
M^{n, f}_t:=f(\omega_t)-
f(\om_0)-\int^t_0(\Delta+b_n\cdot\nabla) f(s,\omega_s)\dif s, \quad t\in [0,T]. 
$$
We want to take weak limits, where the key point is to show
\begin{align}\label{eq-limit1}
\lim_{k\to\infty}\mE^{n_k}_{x}
\left(\int^{t}_0(
b^{n_k}\cdot\nabla f)(s,\omega_s)\dif s\cdot G(\omega)\right)=\mE_{x}
\left(\int^{t}_0(
b\cdot\nabla f)(s,\omega_s)\dif s\cdot G(\omega)\right).
\end{align} 
Assume that supp$(f)\subset B_R$. 
By \eqref{eq-kry1}, we have
\begin{align}\label{eq-limit2}
\begin{split}
&\sup_{n\geq m}\mE^n_{x}\left| \int^{t}_0((b_m-
b_n)\cdot\nabla f)(s,\omega_s)\dif s\cdot G(\omega)\right|\\
\leq& \|G\|_\infty\|\nabla f\|_\infty
\sup_{n\geq m}\mE^n_{x}
\left(\int^{t}_0|(b_m-
b_n)\chi_R^0|(s,\omega_s)\dif s\right)\\
\leq& C \|G\|_\infty\|\nabla f\|_\infty \ \sup_{n\geq m} \|(b_m-
b_n)\chi_R^0\1_{[0,t]}\|_{ \mL^{p_3}_{q_3}}\to 0,\ m\to\infty,  
\end{split}
\end{align} 
where the cutoff function $\chi$ is defined by \eqref{Eq-chi}. Similarly, by \eqref{eq-kry2}, 
\begin{align}\label{eq-limit3}
\mE_{x}\left| \int^{t}_0((b_m-b)\cdot\nabla f)(r,\omega_r)\dif r\cdot G(\omega)\right|\lesssim \|(b_m-
b)\chi_R^0\1_{[0,t]}\|_{\mL^{p_3}_{q_3}}\to 0, \ ( m\to\infty).
\end{align}
On the other hand, for fixed $m\in\mN$, 
$$
\omega\mapsto \int^t_0(b_m\cdot\nabla f)(r,\omega_r)\dif r\cdot G(\omega)\in C_b(C([0,T]; \R^d)),
$$
so we also have
$$
\lim_{k\to\infty}\mE^{n_k}_{x}\left(\int^t_0(b_m\cdot\nabla f)(s,\omega_s)\dif s\cdot G(\omega)\right)
=\mE_{x}\left(\int^t_0(b_m\cdot\nabla f)(s,\omega_s)\dif s\cdot G(\omega)\right),
$$
which together with \eqref{eq-limit2} and  \eqref{eq-limit3} implies \eqref{eq-limit1}. 

{\em Uniqueness:} Let $\mP_x^{(i)}$, $i=1,2$ be two martingale solutions of SDE \eqref{Eq-SDE} and there is a constant $C>0$
such that for all $x\in\mR^d$ and  $f\in\widetilde \mL^p_q(T)$,
\begin{align}\label{eq-krylov}
\mE^{(i)}_x \l(\int_{0}^{T} f(t, \om_t) \d t \r) \leq C \|f\|_{\widetilde \mL^{p_3}_{q_3}(T)}, \quad \forall p_3,q_3\in (1,\infty) \mbox{ with } \frac{d}{p_3}+\frac{2}{q_3}<2. 
\end{align}
Let $(p_3, q_3)$ be the pair of constants in Theorems \ref{Th-PDEmain1} and \ref{Th-PDEmain2} with $d/p_3+2/q_3<2$, respectively.  For any  $f\in C^\infty_c((0,T)\times\mR^d)$, by Theorems \ref{Th-PDEmain1} and \ref{Th-PDEmain2}, there is a unique solution $u\in\widetilde \mH^{2,p_3}_{q_3}(T)$  with $d/p_3+2/q_3<2$ to the following backward equation:
$$
\p_t u+L u+f=0,\ \ u(T)=0,  
$$
where $L:= \Delta + b\cdot\nabla$. Let $u_n(t,x):=u(t,\cdot)*\rho_n(x)$ be the mollifying approximation of $u$. Then we have
$$
\p_t u_n+Lu_n+g_n=0,\ \ u_n(T)=0,
$$
where
$$
g_n=f*\rho_n+(Lu)*\rho_n-L(u*\rho_n).
$$
For $R>0$, define
$$
\tau_R:=\inf\{t\geq 0: |\omega_t|\geq R\}.
$$
By It\^o's formula, we have
\begin{align}\label{JH}
\mE^{(i)} u_n({T\wedge\tau_R}, \omega_{T\wedge\tau_R})=u_n(0,x)-\mE^{(i)} \left(\int^{T\wedge\tau_R}_0g_n(s,\omega_s)\dif s\right),\ \ i=1,2.
\end{align}
From the proofs for Theorems \ref{Th-PDEmain1} and \ref{Th-PDEmain2}, one can see that 
$$
\| (b\cdot \nabla u)-(b\cdot \nabla u)*\rho_n\|_{\widetilde \mL^{p_3}_{q_3}(T)}\to 0, \quad  \|b\cdot \nabla u- b\cdot \nabla u*\rho_n\|_{\widetilde \mL^{p_3}_{q_3}(T)}\to 0. 
$$
Using estimate \eqref{eq-krylov}, we have
\begin{align*}
&\lim_{n\to\infty}\mE^{(i)} \left(\int^{T\wedge\tau_R}_0\Big((Lu)*\rho_n-L(u*\rho_n)\Big)(s,\omega_s)\dif s\right)\\
\leq& C\lim_{n\to\infty} \l\|\chi_R^0\l[ (b\cdot\nabla u)*\rho_n-b\cdot(\nabla u*\rho_n) \r] \r\|_{\mL^{p_3}_{q_3}(T)}\\
\leq & C\lim_{n\to\infty} \l\|\chi_R^0\l[ (b\cdot\nabla u)*\rho_n-b\cdot\nabla u\r]\r\|_{\mL^{p_3}_{q_3}(T)}+C \lim_{n\to\infty}\l\|\chi_R^0\l[b\cdot\nabla u-b\cdot(\nabla u*\rho_n) \r] \r\|_{\mL^{p_3}_{q_3}(T)}=0,
\end{align*}
where the cutoff function $\chi$ is defined by \eqref{Eq-chi}. Recalling that $u\in \tH^{2,p_3}_{q_3}(T)$, $\p_t u\in \tL^{p_3}_{q_3}(T)$ and $d/p_3+2/q_3<2$, due to Lemma 10.2 of \cite{krylov2005strong} $u$ is a bounded H\"older continuous function on $[0,T]\times \R^d$. Letting $n\to\infty$ for both sides of  \eqref{JH} and by the dominated convergence theorem, we obtain
$$
\mE^{(i)} u({T\wedge\tau_R}, \omega_{T\wedge\tau_R})=u(0,x)-\mE^{(i)} \left(\int^{T\wedge\tau_R}_0f(s,\omega_s)\dif s\right),\ \ i=1,2,
$$
which, by letting $R\to\infty$ and noting that $u(T)=0$, yields
$$
u(0,x)=\mE^{(i)} \left(\int^{T}_0f(s,\omega_s)\dif s\right),\ \ i=1,2.
$$
This in particular implies the uniqueness of martingale solutions (see \cite[Corollary 6.2.6]{stroock2007multidimensional}).
\epf

\section{Appendix}
In this section, we present sketches of proofs for Lemmas \ref{Le-energy} and \ref{Le-Max}. 
\begin{proof}[Proof of Lemma \ref{Le-energy}]
As presented in the proof of \cite[Lemma 3.2]{zhao2019stochastic}, for almost every $s, t\in I$ with $s<t$, 
\begin{align}\label{eq-energy}
\begin{aligned}
&\frac{1}{2} \left(\int u_k^2\eta^2 \right)(t)-\frac{1}{2} \left(\int u_k^2\eta^2\right)(s)+  \int_s^t\!\!\!\int \nabla u_k \cdot \nabla (u_k \eta^2)\\
\leq & -\int_s^t\!\!\!\int (u_k+k) b\cdot \nabla (u_k\eta^2) - \int_s^t\!\!\!\int  \div b \, (u_k+k) \, u_k\eta^2+\int_s^t\!\!\!\int fu_k\eta^2.  
\end{aligned}
\end{align}
H\"older's inequality yields  
\be\label{eq-uk'}
\begin{aligned}
\int_s^t\!\!\!\int \nabla u_k \cdot \nabla (u_k \eta^2) =&\int_s^t\!\!\!\int |\nabla u_k \eta|^2 +2  \int_s^t\!\!\!\int (\nabla u_k \eta)\cdot (u_k\nabla \eta )\\
\geq&  \frac{1}{2}\int_s^t\!\!\!\int |\nabla u_k \eta|^2  -\frac{C}{(R-\rho)^2} \|u_k\|_{\mL^2_2(A_s^t(k))}^2
\end{aligned}
\ee

\medskip
Integration by parts and H\"older's inequality yield  
\begin{equation*}
\begin{aligned}
&-\int_s^t\!\!\!\int (u_k+k) b\cdot \nabla (u_k\eta^2)\\
=&- \frac{1}{2}\int_s^t\!\!\!\int \eta^2 b\cdot \nabla( u_k^2) - 2\int_s^t\!\!\!\int u_k^2 \eta b\cdot \nabla \eta- k \int_s^t\!\!\!\int \eta^2 b\cdot \nabla u_k-2 k\int_s^t\!\!\!\int u_k\eta b\cdot \nabla \eta\\
 =& \l[ \int_s^t\!\!\!\int  u_k^2\eta b\cdot \nabla \eta +\frac{1}{2} \int_s^t\!\!\!\int  \div b\  u_k^2\eta^2 \r]-2\int_s^t\!\!\!\int u_k^2 \eta b\cdot \nabla \eta\\
 &+\l[ 2k\int_s^t\!\!\!\int u_k\eta b\cdot\nabla \eta+k \int_s^t\!\!\!\int  \div b\  u_k \eta^2\r] -2 k\int_s^t\!\!\!\int u_k\eta b\cdot \nabla \eta \\
=& - \int_s^t\!\!\!\int  u_k^2\eta b\cdot \nabla \eta +\frac{1}{2} \int_s^t\!\!\!\int  \div b\  u_k^2\eta^2 +k \int_s^t\!\!\!\int  \div b\  u_k \eta^2. 
\end{aligned}
\end{equation*}
Therefore,  
\be
\begin{aligned}
&-\int_s^t\!\!\!\int (u_k+k) b\cdot \nabla (u_k\eta^2) - \int_s^t\!\!\!\int  \div b \, (u_k+k) \, u_k\eta^2 \\
=& - \int_s^t\!\!\!\int  u_k^2\eta b\cdot \nabla \eta -\frac{1}{2} \int_s^t\!\!\!\int  \div b\  u_k^2\eta^2 \\
\leq &\frac{2}{R-\rho} \int_s^t\!\!\!\int \l|  bu_k^2 \eta \r| + \frac12\int_s^t\!\!\!\int  |\div b|\,  u_k^2\eta^2\\
\leq&   \frac{2}{(R-\rho)}  \|b\|_{\mL^{p_2}_{q_2}(Q)}   \|u_k\|_{\mL^{p^*_2}_{q^*_2}(A_s^t(k))}^2+  \|\div b\|_{\mL^{p_2}_{q_2}(Q)}  \|u_k\|_{\mL^{p^*_2}_{q^*_2}(A_s^t(k))}^2. 
\end{aligned}
\ee
By H\"older's inequality, 
\begin{align}\label{eq-fu}
\begin{aligned}
\int_s^t\!\!\!\int fu_k\eta^2 \leq& \|f\|_{\mL^{p_3}_{q_3}(Q)} \|\1_{A_s^t(k)}\|_{\mL^{p^*_3}_{q^*_3}} \|u_k\|_{\mL^{p^*_3}_{q^*_3}(A_s^t(k))}  \\
\leq& \frac{1}{2}\|f\|_{\mL^{p_3}_{q_3}(Q)}^2\|\1_{A_s^t(k)}\|_{\mL^{p^*_3}_{q^*_3}}^2 + \frac{1}{2} \|u_k\|_{\mL^{p^*_3}_{q^*_3}(A_s^t(k))}^2
\end{aligned}
\end{align}
Combing \eqref{eq-energy}-\eqref{eq-fu} and using H\"older's inequality, we obtain \eqref{Eq-energy}. 
\end{proof}

\bl\label{Le-decrease}
Suppose $\{y_j\}_{j\in \N}$ is a nonnegative nondecreasing real sequence, 
$$
y_{j+1}\leq N C^j y_j^{1+\eps} 
$$
with $\eps>0$ and $C>1$. Assume 
$$
y_0\leq N^{-1/\eps } C^{-1/\eps^2}.  
$$
Then $y_j\to 0$ as $j\to \infty$. 
\el

\begin{proof}[Proof of Lemma \ref{Le-Max}]
(i) For any $k\in \mN$, set 
$$
t_k=-\tfrac{1}{2}(1+2^{-k}),\quad B_k'=B_{\frac{1}{2}(1+2^{-k})},\quad  Q'_k=(t_k,0)\times B_k'. 
$$
The cut off functions $\eta_k$ is supported in $B'_{k-1}$ and equals to $1$ in $B'_{k}$ such that $|\nabla^i \eta_k|\leq C 2^{ik}\, (i=0,1,2)$. 
Let $M>0$, which will be determined later, and define 
$$
M_k:=M(2-2^{-k}),\quad u_k:= (u-M_k)^+,\quad U_{k}:= \|u_k\|_{\mL^{2}_2(Q'_k)}^2+\sum_{i=2}^3 \|u_k\|_{\mL^{p^*_i}_{q^*_i}(Q'_k)}^2 
$$
and 
$$
E_k:=\sup_{t\in [t_k, 0]}\int (u_k\eta_k)^2(t)+\int_{t_k}^0\!\!\!\int |\nabla (u_k\eta_k)|^2. 
$$
For any $s,t$ satisfying $t_k\leq s\leq t_{k+1}\leq t\leq 0$, by Lemma \ref{Le-energy}, we have 
\begin{align*}
&\int (u_{k+1}\eta_{k+1})^2 (t) +\int_s^t\!\!\!\int |\nabla(u_{k+1}\eta_{k+1})|^2\\
\leq &\int (u_{k+1}\eta_{k+1})^2 (s)+ C \|f\|^2_{\mL^{p_3}_{q_3}(Q_1)}\|\1_{\{u_{k+1}>0\}\cap Q'_k}\|_{\mL^{p_3^*}_{q_3^*}}^2 \\
& +C^k \left( \|u_{k+1}\|_{\mL^{2}_{2}(Q'_k)}^2+ \sum_{i=2}^3 \|u_{k+1}\|_{\mL^{p^*_i}_{q^*_i}(Q'_k)}^2 \right)
\end{align*}
Using the range of $s, t$ and taking the mean value in $s$ between $t_{k+1}$ and $t_k$, we get 
\begin{equation*} 
\begin{aligned}
&\int (u_{k+1}\eta_{k+1})^2 (t)+\int_{t_{k+1}}^t\!\!\int |\nabla(u_{k+1}\eta_{k+1})|^2\\
\leq & 4\cdot2^{k} \int^{t_{k+1}}_{t_k}\!\!\int (u_{k+1}\eta_{k+1})^2 + C \|f\|^2_{\mL^{p_3}_{q_3}(Q_1)} \|\1_{\{u_{k+1}>0\}\cap Q'_k}\|_{\mL^{p^*_3}_{q^*_3}}^2 \\ 
& +C^k \left( \|u_{k+1}\|_{\mL^{2}_{2}(Q'_k)}^2+ \sum_{i=2}^3\|u_{k+1}\|_{\mL^{p_i^*}_{q_i^*}(Q'_k)}^2 \right)\\
\leq & C^k \left(\|u_{k+1}\|_{\mL^2_2(Q'_k)}^2+\sum_{i=2}^3\|u_{k+1}\|_{\mL^{p^*_i}_{q^*_i}(Q'_k)}^2\right)+ C\|f\|^2_{\mL^{p_3}_{q_3}(Q_1)} \|\1_{\{u_{k+1}>0\}\cap Q'_k}\|_{\mL^{p_3^*}_{q_3^*}}^2. 
\end{aligned}
\end{equation*}
Choosing $M>C\|f\|_{\mL^{p_3}_{q_3}(Q_1)}$, the above inequalities yield, 
\begin{equation}\label{Eq-Ek+1-Uk}
\begin{aligned}
E_{k+1}\leq &  \sup_{t\in [t_{k+1}, 0]} \int (u_{k+1}\eta_{k+1})^2(t)+\int_{t_{k+1}}^0 \!\!\int |\nabla(u_{k+1}\eta_{k+1})|^2\\
\leq & C^k \left(\|u_{k+1}\|_{\mL^2_2(Q'_k)}^2+ \sum_{i=2}^3\|u_{k+1}\|_{\mL^{p^*_i}_{q^*_i}(Q'_k)}^2\right)+ C \|f\|^2_{\mL^{p_3}_{q_3}(Q_1)} \|\1_{\{u_{k+1}>0\}\cap Q'_k}\|_{\mL^{p_3^*}_{q_3^*}}^2 \\
\leq&  C^k U_k+ M^2 \|\1_{\{u_{k+1}>0\}\cap Q'_k}\|_{\mL^{p_3^*}_{q_3^*}}^2. 
\end{aligned}
\end{equation}
The quantity $E_{k+1}$ controls $u_{k+1}\eta_{k+1}$ in $\mL_\infty^{2}(Q'_{k+1})$,  and thanks to Sobolev embedding, also in the space $\mL_2^{\gamma}(Q'_{k+1})$ for $\gamma=2d/(d-2)$ if $d\geq 3$ and $\mL_2^{\gamma}(Q'_{k+1})$ for any $\gamma\in [2,\infty)$ if $d=2$. 
So, by interpolation, $E_{k+1}$ controls the $\mL^r_s(Q'_{k+1})$-norm of $u_{k+1}\eta_{k+1}$ for any 
\be\label{eq-Lrs}
r,s \geq 2 \ \mbox{ with }\ \frac{d}{r}+\frac{2}{s}>\frac{d}{2}.  
\ee 
Noting that $d/p_i+2/q_i<2\ (i=2,3)$   and \eqref{Eq-*-pq}, we have 
$$
\frac{d}{p_i^*}+\frac{2}{q_i^*}>\frac{d}{2}, \ (i=2,3). 
$$
By H\"older's inequality and \eqref{Eq-Ek+1-Uk}, one can see that there exists a constant $\eps>0$ such that 
\begin{align*}
U_{k+1}= & \|u_{k+1}\|_{\mL^{2}_2(Q'_{k+1})}^2+\sum_{i=2}^3\|u_{k+1}\|_{\mL^{p_i^*}_{q_i^*}(Q'_{k+1})}^2\\
\leq &C E_{k+1} |\{u_{k+1}>0\}\cap Q'_{k+1}|^{2\eps}\\
\leq & C^k U_k |\{u_{k+1}>0\}\cap Q'_k|^{2\eps}+ M^2  \|\1_{\{u_{k+1}>0\}\cap Q'_k}\|_{\mL^{p^*_3}_{q^*_3}}^2 |\{u_{k+1}>0\}\cap Q'_k|^{2\eps}.
\end{align*}
On the other hand, 
\begin{align*}
& M^2  \|\1_{\{u_{k+1}>0\}\cap Q'_k}\|_{\mL^{p^*_3/2}_{q^*_3/2}} \leq  M^2 (M2^{-k-1})^{-1}\|u_k\1_{\{u_k>2^{-k-1}M\}\cap Q'_k}\|_{\mL^{p^*_3}_{q^*_3}}^2\\
\leq& 2^{k+1}M \|u_k\|_{\mL^{p^*_3}_{q^*_3}(Q_k')} \|\1_{\{u_k>2^{-k-1}M\}\cap Q'_k}\|_{\mL^{p^*_3}_{q^*_3}}^2\\
\leq& 4\cdot 4^k \|u_k\|^2_{\mL^{p^*_3}_{q^*_3}(Q_k')}\leq C^k U_k, 
\end{align*}
hence, 
\const{\CUk}
\begin{align*}
U_{k+1}\leq & C^k U_k |\{u_{k+1}>0\}\cap Q'_k|^{2\eps} \leq C^k U_k \left[2^{k+1} M^{-1} \|u_k\|_{\mL^2_2(Q'_k)} \right]^{2\eps}\\
\leq  & M^{-2\eps} C_{\CUk}^k U_k^{1+\eps}.  
\end{align*}
Choosing 
$$
M:=C_{\CUk}\|f\|_{\mL^{p_3}_{q_3}(Q_1)}+C_{\CUk}^{1/(2\eps^2)}\big(\|u^+\|_{\mL^2_2(Q_1)}+\sum_{i=2}^3 \|u^+\|_{\mL^{p_i^*}_{q_i^*}(Q_1)}\big), 
$$ 
we have 
$$
U_0\leq \l(\|u^+\|_{\mL^2_2(Q_1)}+\sum_{i=2}^3\|u^+\|_{\mL^{p_i^*}_{q_i^*}(Q_1)}\r)^2 \leq M^2 C_{\CUk}^{-1/\eps^2}. 
$$
By Lemma \ref{Le-decrease}, 
$$
\|(u-2M)^+\|_{\mL^2_2(Q_{1/2})}\leq \lim_{k\to \infty} U_k =0.  
$$
By the definition of $M$, we obtain 
$$
\|u^+\|_{L^\infty(Q_{1/2})} \leq 2M \leq C \left(\|u^+\|_{\mL^2_2(Q_1)}+\sum_{i=2}^3\|u^+\|_{\mL^{p_i^*}_{q_i^*}(Q_1)}+\|f\|_{\mL_{q_3}^{p_3}(Q_1)}\right). 
$$
\end{proof}

\const{\CGlobalMax}
\const{\Cglobalmax}
\bl\label{Le-GlobalMax}
Let $d\geq 3$, $p_i, q_i\in (1,\infty)$, $i=1,2,3$ and $f\in  \tL^{p_3}_{q_3}(T)$ with $d/p_3+2/q_3<2$. 
Assume $b$ satisfies one of the following two conditions:  
\begin{enumerate}[(a)]
\item   $b=b_0+b_1$,  $\|b_0\|_{\tL^{d,\infty}_\infty}\leq \eps(d)$, for some $\eps(d)>0$ only depending on $d$,  and $b_1\in  \mL^{p_1}_{q_1}$ with $d/p_1+2/q_1=1$ and $p_1\in (d, \infty] $, 
\item  $b, \div b\in  \tL^{p_2}_{q_2}$ with $p_2, q_2\in [2,\infty)$ and $d/p_2+2/q_2<2$.
\end{enumerate}
Then equation \eqref{Eq-PDE} admits a unique weak solution $u\in \widetilde V^0(T)\cap \mL^\infty(T)$. Moreover, the following estimate is valid:
\begin{align}\label{Eq-Max}
\|u\|_{\widetilde V(T)}+\|u\|_{\mL^\infty(T)} \leq C_{\CGlobalMax}  \|f\|_{\widetilde \mL^{p_3}_{q_3}(T)}, \tag{\bf{GM}}
\end{align}
where the constant $C_{\CGlobalMax}$ only depends on  $d, p_1, q_1, p_3, q_3, \eps, T$ and $b_1$ for the first case and $d, p_2, q_2$,  $p_3, q_3, T, \|b\|_{\tL^{p_2}_{q_2}}$ and $\|\div b\|_{\tL^{p_2}_{q_2}}$ for the second case. 
\el
\bpf
Here we only give the proof for the first case, since the second case was essentially proved in \cite{zhang2020stochastic} and \cite{zhao2019stochastic}.  

For any $x\in \R^d$, let $\eta\in C_c^\infty(B_1)$ such that $\eta\equiv1$ on $B_{\frac{1}{2}}$ and $\eta_x:= \eta(\cdot-x)$. As presented in the proof of \cite[Lemma 3.2]{zhao2019stochastic}, for almost every $t\in [0,T]$, 
\begin{align}\label{eq-energy2}
\begin{aligned}
&\frac{1}{2} \left(\int  u_k^2\eta_x^2 \right)(t)+ \int_0^t\!\!\!\int  \nabla u_k \cdot \nabla (u_k \eta_x^2)\leq  \int_0^t\!\!\!\int   b\cdot \nabla u\,   (u_k\eta_x^2) +\int_0^t\!\!\!\int  fu_k\eta_x^2.  
\end{aligned}
\end{align}
As showed in \eqref{eq-uk'} and \eqref{eq-fu}, we have 
\be
\int_0^t\!\!\!\int  \nabla u_k \cdot \nabla (u_k \eta_x^2) \geq  \frac{1}{2} \|\nabla u_k \eta_x\|^2_{\mL^2_2(t)}  -C \|u_k\|_{\mL^2_2(Q(t,x))}^2, 
\ee
and 
\be
\begin{aligned}
\int_0^t\!\!\!\int  fu_k\eta_x^2 \leq& \|f\|_{\tL^{p_3}_{q_3}(T)} \|\1_{A(t,x;k)}\|_{\mL^{p^*_3}_{q^*_3}} \|u_k\eta_x\|_{\mL^{p^*_3}_{q^*_3}(Q(t,x))} \\
\leq & \frac{1}{15} \|u_k\eta_x\|_{V(t)}^2 +C  \|f\|_{\tL^{p_3}_{q_3}(T)}^2 \|\1_{A(t,x;k)}\|_{\mL^{p^*_3}_{q^*_3}}^2, 
\end{aligned}
\ee
where $Q(t,x)= (0,t)\times B_1(x)$ and $A(t,x;k)=\{u>k\}\cap Q(t,x)$. 
Let $b^N_1= (-N\vee b_1)\wedge N$. Then $ \delta_N:=\|b-b^N_1\|_{\mL^{p_1}_{q_1}(T)}\to 0, \ (N\to \infty)$. Furthermore, 
\be
\begin{aligned}
\int_0^t\!\!\!\int   b\cdot \nabla u\,   (u_k\eta_x^2)  =& \int_0^t\!\!\!\int   b_0 \cdot (\nabla u_k\eta_x)\,   (u_k\eta_x) + \int_0^t\!\!\!\int    (b-b_1^N) \cdot (\nabla u_k\eta_x)\,   (u_k\eta_x) \\
&+ \int_0^t\!\!\!\int  b_1^N \cdot (\nabla u_k\eta_x)\,   (u_k\eta_x) =: I_1+ I_2+ I_3.  
\end{aligned}
\ee
For $I_1$,  by \cite[Exercise 1.4.19]{grafakos2008classical} and  \cite[Remark 5]{tartar1998imbedding}, we have 
\const{\CLorentz}
\begin{align*}
I_1\leq& \|b_0\|_{\tL^{d,\infty}_\infty(T)}  \|\nabla u_k \eta_x\|_{\mL^2_2(t)} \|u_k\eta_x\|_{L^2([0,t];\, L^{\frac{2d}{d-2} ,2})}\\
\leq& C_{\CLorentz} \eps \l(\|\nabla u_k \eta_x\|_{\mL^2_2(t)}^2+  \| u_k \nabla \eta_x\|_{\mL^2_2(t)}^2\r).  
\end{align*}
Here $L^{p,q}$ is the Lorentz space and   $C_{\CLorentz}=C_{\CLorentz}(d)$. Choosing $\eps=\eps(d)>0$ small, then 
\be
I_1\leq \frac{1}{15} \|\nabla u_k \eta_x\|_{\mL^2_2(t)}^2 +C \|u_k\|_{\mL^2_2(Q(t,x))}^2. 
\ee
For $I_2$, we have 
$$
I_2\leq  \|b-b^N_1\|_{\mL^{p_1}_{q_1}(T)}  \|\nabla u_k \eta_x\|_{\mL^2_2(t)}  \|u_k \eta_x\|_{\mL^{p_1'}_{q_1'}(t)}, 
$$
where $\frac{1}{p'_1}=\frac{1}{2}-\frac{1}{p_1}$ and $\frac{1}{q'_1}=\frac{1}{2}-\frac{1}{q_1}$. Noting that $d/p_1'+2/q_1'=d/2$,  by choosing $N$ sufficiently large, we get  
\be
I_2\leq C\delta_N \l( \|\nabla u_k \eta_x\|_{\mL^2_2(t)} \|u_k\eta_x\|_{V(t)}\r)\leq \frac{1}{15}\|\nabla u_k \eta_x\|_{\mL^2_2(t)}^2+ \frac{1}{15}\|u_k\eta_x\|_{V(t)}^2. 
\ee
For $I_3$, by H\"older's inequality, we also have 
\be\label{eq-bI3}
I_3\leq N \|\nabla u_k \eta_x\|_{\mL^2_2(t)} \| u_k \eta_x\|_{\mL^2_2(t)} \leq \frac{1}{15} \|\nabla u_k \eta_x\|_{\mL^2_2(t)}^2+ C_N \|u_k\|_{\mL^2_2(Q(t,x))}^2
\ee
Combing \eqref{eq-energy2}-\eqref{eq-bI3}, we obtain 
\be\label{eq-energy-uk}
\|u_k\eta_x\|_{V(t)}\leq C \l( \|u_k\|_{\mL^2_2(Q(t,x))} + \|f\|_{\tL^{p_3}_{q_3}(T)} \|\1_{A(t,x;k)}\|_{\mL^{p^*_3}_{q^*_3}} \r), 
\ee
where $C$ only depends on $d, p_i ,q_i, \eps$ and $b_1$. Now our desired results can be obtained in the same way as in the proofs for Theorem 3.4 and Theorem 3.6 in \cite{zhao2019stochastic}. 
\epf

\section*{Acknowledgement}
The second named author is very grateful to Professor Nicolai Krylov and Xicheng Zhang who encouraged him to persist in studying this problem, and also Professor Kassmann for providing him an excellent environment to work at Bielefeld University. 

\bibliographystyle{alpha2}

\bibliography{mybib}
\end{document}